\documentclass[oneside]{amsart}
\usepackage{amssymb, a4wide, mathdots, url,hyperref,graphicx}

\usepackage{amsfonts}
\usepackage{amssymb}
\usepackage{amsfonts, amsmath, amsthm, amssymb}
\usepackage{graphicx}
\usepackage{listings}
\usepackage{xcolor}
\usepackage{thmtools}
\usepackage{mathrsfs}
\usepackage{hyperref}
\usepackage{marginnote}
\usepackage[all]{xy}
\usepackage{accents}
\usepackage{enumerate}

\newcommand{\diag}{\operatorname{diag}}

\newcommand{\rk}{\operatorname{rk}}

\theoremstyle{plain}
\newtheorem{thm}{Theorem}[section]

\newtheorem{prop}[thm]{Proposition}
\newtheorem{lem}[thm]{Lemma}
\newtheorem{cor}[thm]{Corollary}

\newtheorem{defn}[thm]{Definition}
\newtheorem{defn/thm}[thm]{Definition/Theorem}

\theoremstyle{remark}

\include{psfig}
\begin{document}
	
	\title{Irreducibility of Polarized Automorphic Galois Representations of Type $G_2$}
	
	\author{Boyi Dai}
    \address{Beijing International Center for Mathematical Research, Peking University, Beijing, China}
	\email{daiboymath@gmail.com}
	
	\begin{abstract}
		Let $\pi$ be a regular algebraic polarized cuspidal automorphic representation of $\mathrm{GL}_7$ over a totally real or CM field, and let $\{\rho_{\lambda}\}_{\lambda}$ be the associated compatible system. We prove that if for some coefficient place $\lambda_0$, the derived monodromy group of $\rho_{\lambda_0}$ is of type $G_2$, then there exists a Dirichlet density one set $\mathcal{L}$ of rational primes such that $\rho_{\lambda}$ is irreducible for every $\ell\in\mathcal{L}$ and every $\lambda|\ell$. Combining this result with the work of Feng and Whitmore and treating a few additional cases, we remove the numerical restriction $7\nmid n$ in their density-one irreducibility theorem for polarized automorphic compatible systems. As further applications, we give a density-one large residual image result for automorphic compatible systems of type $G_2$ and realize $G_2(\mathbb{F}_{\ell})^m$ for $1\le m\le11$ as Galois groups over $\mathbb{Q}$ unramified outside $\ell$, for a density-one set of primes $\ell$.
	\end{abstract}
	
	\maketitle
	

	\section{Introduction}
	
	A fundamental prediction of the Langlands program is that the Galois representations associated with cuspidal automorphic representations of general linear groups should be irreducible. Let $K$ be a CM field or a totally real field, and let $(\pi,\chi)$ be a regular algebraic polarized cuspidal automorphic representation of $\mathrm{GL}_n(\mathbb{A}_K)$. By the work of many authors (see \cite[Theorem 2.1.1]{BLGGT14}), one can attach to $\pi$ an $E$-rational compatible system of semisimple Galois representations
	$$\{\rho_{\lambda}:\mathrm{Gal}_K\to\mathrm{GL}_n(\overline{E}_{\lambda})\}_{\lambda}.$$
	It is conjectured that $\rho_{\lambda}$ is irreducible for every $\lambda$.
	
	Some of the earliest results towards this conjecture are due to Ribet \cite{Ri77}, Taylor \cite{Ta95}, and Blasius and Rogawski \cite{BR92}, who treated the two- and three-dimensional cases of totally real fields in the polarized setting. Calegari and Gee \cite{CG13} established density-one irreducibility for $n\leq5$ of totally real fields. In arbitrary dimension, Barnet-Lamb, Gee, Geraghty and Taylor \cite{BLGGT14} proved a density-one irreducibility result of CM fields or totally real fields, under an additional assumption on Hodge-Tate weights. Without this additional assumption, Patrikis and Taylor \cite{PT15} obtained a positive-density result on arbitrary dimension of CM fields or totally real fields. Xia \cite{Xi19} proved the density-one result to $n\leq6$ of imaginary CM fields, and Hui \cite{Hu23b} strengthened this to irreducibility at all but finitely many coefficient places of CM fields or totally real fields. B\"ockle and Hui \cite{BH25} established the conjecture for $n=3$ over totally real fields without assuming polarization. Hui and Lee \cite{HL24,HL26} obtained irreducibility results in arbitrary dimension under the connectedness assumption on the compatible systems, and certain type $A_1$ hypotheses on the algebraic monodromy groups. More recently, Feng and Whitmore \cite{FW25} proved density-one irreducibility results in the polarized setting under the numerical assumption that $7\nmid n$ and that, if $4\mid n$, then $n=4p$ for some prime $p$. There is also a variant of this conjecture in the setting of strictly compatible systems; see \cite{DW21,Hu23a,Dai25,DWW26} for results in low dimensions.
	
	The case of seven-dimensional representation with derived monodromy group of type $G_2$ remains open in the preceding work. The difficulty in this case arises from the fact that the restriction of its seven-dimensional irreducible representation to the equal-rank subgroup $\mathrm{SL}_3$ decomposes as $\mathrm{Std}\oplus\mathrm{Std}^{\vee}\oplus\mathbf{1}$. Since the complex conjugation may fix the two three-dimensional constituents, the derived formal character together with the information on complex conjugation does not suffice to rule out this possibility. In this paper, we give a density-one irreducibility result in this case.
	
	\begin{thm}\label{main}
		Let $K$ be imaginary CM or totally real, and let $\{\rho_{\lambda}:\mathrm{Gal}_K\to\mathrm{GL}_n(\overline{E}_{\lambda})\}_{\lambda}$ be the $E$-rational compatible system attached to a regular algebraic polarized cuspidal automorphic representation $(\pi,\chi)$ of $\mathrm{GL}_n(\mathbb{A}_K)$. Suppose that for some $\lambda_0$, the representation $\rho_{\lambda_0}$ is Lie-irreducible, and the Lie algebra of the derived monodromy group $\mathbf{G}_{\lambda_0}^{\mathrm{der}}$, together with its tautological representation, is isomorphic to $\mathfrak{g}_2$ with its standard seven-dimensional irreducible representation. Then there exists a set $\mathcal{L}$ of rational primes of Dirichlet density one such that $\rho_{\lambda}$ is irreducible for every $\ell\in\mathcal{L}$ and every $\lambda|\ell$.
	\end{thm}
	
	Our proof is inspired by Larsen and Pink \cite[Proposition 8.2 and 8.8]{LP92}. In their work, they use an openness result of Serre to construct arbitrarily many Frobenius tori whose splitting fields are linearly disjoint over a fixed number field, and then apply the Chebotarev density theorem to obtain a density-one result. In our setting, we use the openness theorem of Pink \cite{Pi98} to construct arbitrarily many linearly disjoint quadratic extensions and apply the Chebotarev density theorem to show that the bad primes have Dirichlet density zero.
	
	Combining our main theorem with the work of \cite{FW25}, together with a few cases not explicitly treated in their paper but covered by their methods, we obtain the following extension of their main result, removing the assumption that $7\nmid n$.
	\begin{thm}\label{main2}
		Let $K$ be totally real or imaginary CM. Let $(\pi,\chi)$ be a polarized regular algebraic cuspidal automorphic representation of $\mathrm{GL}_n(\mathbb{A}_K)$ and $\{\rho_{\lambda}\}_{\lambda}$ be its associated compatible system of Galois representations. Suppose that if $4|n$, then $n=4p$ for some prime number $p$. Then there exists a Dirichlet density one set of rational primes $\mathcal{L}$ such that if $\ell\in\mathcal{L}$ and $\lambda|\ell$, then $\rho_{\lambda}$ is irreducible.
	\end{thm}
	
	As a further application of our main theorem, Larsen's maximality theorem \cite[Theorem 3.17]{La95} yields the following large image result for the residual representations of automorphic compatible systems of type $G_2$. This extends the results of Zenteno \cite{Zen23}.
	\begin{cor}\label{Argon}
		Let $K$ be imaginary CM or totally real, and let $\pi$ be a polarized regular algebraic cuspidal automorphic representation of $\mathrm{GL}_7(\mathbb{A}_K)$ with associated $E$-rational compatible system
		$$\{\rho_{\lambda}:\mathrm{Gal}_K\to\mathrm{GL}_7(\overline{E}_{\lambda})\}_{\lambda}.$$
		Suppose that for some coefficient prime $\lambda_0$, the derived monodromy group of $\rho_{\lambda_0}$ is of type $G_2$. Then there exists a set $\mathcal{L}$ of rational primes of Dirichlet density one such that, for every $\ell\in\mathcal{L}$ and every $\lambda|\ell$, the projective image of the semisimple reduction $\overline{\rho}_{\lambda}$ is isomorphic to $G_2(\mathbb{F}_{\ell^{s_{\lambda}}})$ for some positive integer $s_{\lambda}$.
	\end{cor}
	
	The above corollary, together with Chenevier's construction of level-one cuspidal automorphic representations of type $G_2$ in \cite{Ch19}, yields the following density-one result on the inverse Galois problem.
	\begin{cor}\label{Calcium}
		There exists a Dirichlet density one set $\mathcal{L}$ of rational primes such that, for every $\ell\in\mathcal{L}$ and every $1\le m\le 11$, there is a Galois extension $L/\mathbb{Q}$, unramified at every rational prime other than $\ell$, with
		$$\mathrm{Gal}(L/\mathbb{Q})\simeq G_2(\mathbb{F}_{\ell})^m.$$
	\end{cor}
	
	Feit and Fong \cite{FF85} and Thompson \cite{Th85} proved regular realizations of $G_2(\mathbb{F}_{\ell})$ for every prime $\ell\ge5$. The question of realizing such groups with ramification at only one prime has also been studied. In particular, Lansky and Pollack \cite{LP02} predicted the existence of a Galois extension of $\mathbb{Q}$ with Galois group $G_2(\mathbb{F}_5)$ ramified only at 5; see also \cite{DGV11}. In type $A_1$, Wiese \cite{Wie13} obtained density-one results for $\mathrm{PSL}_2(\mathbb{F}_{\ell^d})$ and $\mathrm{PGL}_2(\mathbb{F}_{\ell^d})$, for even and odd $d$, respectively, assuming Maeda's conjecture. Our corollary gives an unconditional density-one result of this kind for $G_2(\mathbb{F}_{\ell})$ and its direct powers up to the eleventh power.
	
	The article is organized as follows. In Section 2, we recall some preliminary results on compatible systems, including $\ell$-independence properties, algebraic envelopes, and polarization. In Section 3, we prove the main theorem. Section 4 is devoted to the extension of Feng--Whitmore's irreducibility theorem. Finally, in Section 5, we prove the large residual image result and its application to the inverse Galois problem.

	\section{Preliminaries}
	
	\subsection{Compatible systems and $\lambda$-independence}
	
	\begin{defn}
		Let $K$ be a number field and let $F$ be either an $\ell$-adic field or the algebraic closure $\overline{\mathbb{Q}}_{\ell}$. Let
		$$\rho:\mathrm{Gal}_K\to\mathrm{GL}_n(F)$$
		be a semisimple $\ell$-adic Galois representation of $K$. We define the following algebraic groups:
		\begin{itemize}
			\item Algebraic monodromy group $\mathbf{G}$: the Zariski closure of $\rho(\mathrm{Gal}_K)$ in $\mathrm{GL}_{n,F}$.
			\item The connected monodromy group $\mathbf{G}^{\circ}$: the identity component of $\mathbf{G}$.
			\item The derived monodromy group $\mathbf{G}^{\mathrm{der}}$: the derived subgroup of $\mathbf{G}^{\circ}$.
		\end{itemize}
	\end{defn}
	
	Since $\rho$ is semisimple, $\mathbf{G}^{\circ}$ is reductive, and the derived subgroup $\mathbf{G}^{\mathrm{der}}$ is semisimple. For a number field $K$, denote by $\Sigma_K$ the set of its finite places.
	\begin{defn}
	Let $K$ be a number field. An $n$-dimensional $E$-rational Serre compatible system of $\mathrm{Gal}_K$ consists of the data
	$$\mathcal{M}=\left(E,S,\{P_v(T)\}_{v\not\in S},\{\rho_{\lambda}\}_{\lambda\in\Sigma_E}\right),$$
	where:
	\begin{itemize}
		\item $E$ is a number field;
		\item $S$ is a finite set of finite places of $K$, called the exceptional set;
		\item for each finite place $v$ of $K$ outside $S$, $P_v(T)\in E[T]$ is a monic polynomial of degree $n$;
		\item for each $\lambda\in\Sigma_E$, $\rho_{\lambda}:\mathrm{Gal}_K\to\mathrm{GL}_n(\overline{E}_{\lambda})$ is a continuous semisimple $\lambda$-adic representation;
	\end{itemize}
	
	and the following conditions hold:
	\begin{enumerate}[(i)]
		\item For each $\lambda\in\Sigma_E$, $\rho_{\lambda}$ is unramified outside $S\cup S_{\lambda}$, where $S_{\lambda}$ is the set of finite places of $K$ above the residue characteristic of $\lambda$.
		\item For each $\lambda\in\Sigma_E$ and each $v\not\in S\cup S_{\lambda}$, the characteristic polynomial of $\rho_{\lambda}(\mathrm{Frob}_v)$ is $P_v(T)$. 
	\end{enumerate}
    \end{defn}
    
    Thanks to the work of many authors, including, among others, Bella\"iche, Caraiani, Chenevier, Clozel, Harris, Kottwitz, Labesse, Shin, and Taylor, compatible systems of Galois representations can be attached to certain cuspidal automorphic representations.
    \begin{thm}\cite[Theorem 2.1.1]{BLGGT14}\label{Titanium}
    	Let $K$ be a totally real or CM number field, and let $(\pi,\chi)$ be a regular algebraic polarized cuspidal automorphic representation of $\mathrm{GL}_n(\mathbb{A}_K)$ of weight $a$. Then there exist an integer $w$, a number field $E$, and an $n$-dimensional $E$-rational Serre compatible system
    	$$\{\rho_{\pi,\lambda}:\mathrm{Gal}_K\to\mathrm{GL}_n(\overline{E}_{\lambda})\}_{\lambda\in\Sigma_E}$$
    	satisfying the following properties. Write $\ell$ for the residue characteristic of $\lambda$, and fix an isomorphism $\iota_{\lambda}:\overline{E}_{\lambda}\simeq\mathbb{C}$ compatible with the chosen embedding of $E$ into $\mathbb C$.
    	\begin{enumerate}[(i)]
    		\item The pair $(\rho_{\pi,\lambda},\varepsilon_{\ell}^{1-n}\rho_{\chi,\lambda})$ is polarized and totally odd.
    		\item At each finite place $v\nmid\ell$,
    		$$\iota_{\lambda}\mathrm{WD}(\rho_{\pi,\lambda}|_{\mathrm{Gal}_{K_v}})^{\mathrm{F\text{-}ss}}\simeq\mathrm{rec}_{K_v}(\pi_v\otimes|\det|_v^{(1-n)/2}),$$
    		and the corresponding Weil--Deligne representation is pure of weight $w$.
    		\item The representation $\rho_{\pi,\lambda}$ is de Rham at all places above $\ell$. For each embedding $\tau:K\hookrightarrow\overline{E}_{\lambda}$, its Hodge--Tate weights are
    		$$\mathrm{HT}_{\tau}(\rho_{\pi,\lambda})=\{a_{\tau,i}+n-i:1\le i\le n\},$$
    		with
    		$$\mathrm{HT}_{\tau\circ c}(\rho_{\pi,\lambda})=\{w-h:h\in\mathrm{HT}_{\tau}(\rho_{\pi,\lambda})\}.$$
    		\item If $v\mid\ell$ and $\pi_v$ admits an Iwahori-fixed vector, the local compatibility formula in (ii) also holds at $v$. In this case $\rho_{\pi,\lambda}$ is semistable at $v$; it is crystalline when $\pi_v$ is unramified.
    	\end{enumerate}
    \end{thm}
    
    To describe $\lambda$-independent properties of compatible systems, we recall the following notions.
    \begin{defn}
    	Let $F$ be a field and let $\mathbf{G}\subseteq\mathrm{GL}_{n,F}$ be a connected reductive subgroup.
    	\begin{enumerate}[(i)]
    		\item Let $\textbf{T}$ be a maximal torus of $\mathbf{G}\times_F\overline{F}$. The formal character of $\mathbf{G}$ is the conjugacy class of $\textbf{T}$ in $\mathrm{GL}_{n,\overline{F}}$.
    		\item Let $\textbf{T}'$ be a maximal torus of $\mathbf{G}^{\mathrm{der}}\times_F\overline{F}$ contained in $\textbf{T}$. The formal bi-character of $\mathbf{G}$ is the conjugacy class of the chain $\textbf{T}'\subseteq\textbf{T}$ in $\mathrm{GL}_{n,\overline{F}}$.
    		\item Let $F_1,F_2$ be two fields, and let $\mathbf{G}_i\subseteq\mathrm{GL}_{n_i,F_i}, i=1,2$ be two reductive groups. We say that they have the same formal character (resp. formal bi-character) if $n_1=n_2=n$ and there exists a split $\mathbb{Z}$-subtorus $\textbf{T}_{\mathbb{Z}}\subseteq\mathrm{GL}_{n,\mathbb{Z}}$ (resp. a chain of split $\mathbb{Z}$-subtori $\textbf{T}_{\mathbb{Z}}'\subseteq\textbf{T}_{\mathbb{Z}}\subseteq\mathrm{GL}_{n,\mathbb{Z}}$) such that $\textbf{T}_{\mathbb{Z}}\times_{\mathbb{Z}}\overline{F}_i$ (resp. $\textbf{T}_{\mathbb{Z}}'\times_{\mathbb{Z}}\overline{F}_i\subseteq\textbf{T}_{\mathbb{Z}}\times_{\mathbb{Z}}\overline{F}_i$) represents the formal character (resp. formal bi-character) of $\mathbf{G}_i$ for each $i$. This defines an equivalence relation on formal characters (resp. formal bi-characters) of reductive subgroups of general linear groups over arbitrary fields.
    		\item Let $\{F_i\}$ be a family of fields, and let $\{\mathbf{G}_i\subseteq\mathrm{GL}_{n,F_i}\}$ be a family of reductive groups. We say that they have the same formal character (resp. same formal bi-character) if their formal characters (resp. formal bi-characters) belong to a single equivalence class under the relation defined in (iii). We say that they have bounded formal characters (resp. bounded formal bi-characters) if their formal characters (resp. formal bi-characters) belong to finitely many such equivalence classes.
    	\end{enumerate}
    \end{defn}
    
    We call the formal characters of $\mathbf{G}_{\lambda}^{\circ}$ and $\mathbf{G}_{\lambda}^{\mathrm{der}}$ the formal character and the derived formal character of $\rho_{\lambda}$, respectively. We recall the following fundamental $\lambda$-independence results on algebraic monodromy groups.
    \begin{thm}\cite{Se81}, \cite{Se84}, \cite[Theorem 3.19, Remark 3.22]{Hu13}\label{Hydrogen}
    	Let $\{\rho_{\lambda}:\mathrm{Gal}_K\to\mathrm{GL}_n(E_{\lambda})\}_{\lambda}$ be an $E$-rational Serre compatible system. Denote by $\mathbf{G}_{\lambda}$ the algebraic monodromy group of $\rho_{\lambda}\otimes\overline{E}_{\lambda}$.
    	\begin{enumerate}[(i)]
    		\item The kernel of the natural surjection
    		$$\mathrm{Gal}_K\to\pi_0(\mathbf{G}_{\lambda})=\mathbf{G}_{\lambda}/\mathbf{G}_{\lambda}^{\circ}$$
    		is independent of $\lambda$. We call the finite Galois extension of $K$ corresponding to this common kernel the \emph{connectedness field} of the compatible system. In particular, the connectedness of $\mathbf{G}_{\lambda}$ is independent of $\lambda$. 
    		\item The formal bi-character of the tautological representation $\mathbf{G}_{\lambda}^{\circ}\hookrightarrow\mathrm{GL}_{n,\overline{E}_{\lambda}}$, and hence the formal character and the derived formal character of $\rho_{\lambda}$, are independent of $\lambda$. In particular, the rank and semisimple rank of $\mathbf{G}_{\lambda}$ are independent of $\lambda$.
    	\end{enumerate}
    \end{thm}
    
    We recall the following recent result of B\"ockle and Hui on the descent of coefficients of compatible systems.
    \begin{thm}\cite[Theorem 1.3(i) and (iii)]{BH26}\label{Scandium}.
    	Let
    	$$\{\rho_{\lambda}:\mathrm{Gal}_K\to\mathrm{GL}_n(\overline{E}_{\lambda})\}_{\lambda}$$
    	be an $E$-rational Serre compatible system of a number field $K$. After replacing $E$ by a suitable finite extension, we may assume that, after a change of coordinates, each $\rho_{\lambda}$ is defined over $E_{\lambda}$ and every irreducible constituent is absolutely irreducible.
    \end{thm}
	
	\subsection{Algebraic envelope}
	
	Denote by $(\overline{\rho}^{\mathrm{ss}},\overline{V}^{\mathrm{ss}})$ the semisimple reduction of a $\lambda$-adic Galois representation $(\rho,V)$. Let $\{(\rho_{\lambda},V_{\lambda})\}_{\lambda}$ be an $n$-dimensional $E$-rational Serre compatible system of a number field $K$. By \autoref{Scandium}, after enlarging $E$ by a suitable finite extension, we may write the system as $\{\rho_{\lambda}:\mathrm{Gal}_K\to\mathrm{GL}_n(E_{\lambda})\}_{\lambda}$. Put $d=[E:\mathbb{Q}]$. By restriction of scalars, we obtain an $nd$-dimensional $\mathbb{Q}$-rational compatible system:
	$$\left\{\rho_{\ell}:=\bigoplus_{\lambda|\ell}\rho_{\lambda}:\mathrm{Gal}_K\to\left(\mathrm{Res}_{E/\mathbb{Q}}\mathrm{GL}_{n,E}\right)\left(\mathbb{Q}_{\ell}\right)\subseteq\mathrm{GL}_{nd}(\mathbb{Q}_{\ell})\right\}_{\ell}.$$
	\begin{defn/thm}\cite[Theorem 3.1]{Hu23b}\label{Vanadium}
		Suppose that there exist integers $N_1,N_2\ge0$ and a finite extension $K'/K$ such that the following conditions hold.
		\begin{enumerate}
		 	\item[(a)] (Bounded tame inertia weights): For all $\ell\gg0$ and each finite place $v$ of $K$ above $\ell$, the tame inertia weights of the local representation $(\overline{\rho}_{\ell}^{\mathrm{ss}}\otimes\overline{\varepsilon}_{\ell}^{N_1})|_{\mathrm{Gal}_{K_v}}$ belong to $[0,N_2]$.
		 	\item[(b)] (Potential semistability): For all $\ell\gg0$ and each finite place $w$ of $K'$ not above $\ell$, the semisimplification of the local representation $\overline{\rho}_{\ell}^{\mathrm{ss}}|_{\mathrm{Gal}_{K'_w}}$ is unramified.
		\end{enumerate}
		Then there exists a finite Galois extension $L/K$ and for each $\ell\gg0$, a connected reductive subgroup $\underline{G}_{\ell}\subseteq\mathrm{GL}_{nd,\mathbb{F}_{\ell}}$ satisfying
		\begin{enumerate}[(i)]
			\item The derived group $\underline{G}_{\ell}^{\mathrm{der}}$ is the Nori group of $\overline{\rho}_{\ell}^{\mathrm{ss}}(\mathrm{Gal}_K)\subseteq\mathrm{GL}_{nd}(\mathbb{F}_{\ell})$.
			\item The image $\overline{\rho}_{\ell}^{\mathrm{ss}}(\mathrm{Gal}_L)$ is a subgroup of $\underline{G}_{\ell}(\mathbb{F}_{\ell})$ whose index is uniformly bounded as $\ell$ varies.
			\item $\underline{G}_{\ell}$ acts semisimply on the ambient space.
			\item The formal characters of the embeddings $\underline{G}_{\ell}\hookrightarrow\mathrm{GL}_{nd,\mathbb{F}_{\ell}}$ are bounded as $\ell$ varies.
		\end{enumerate}
		For all but finitely many $\lambda$, let $\ell$ be the rational prime below $\lambda$ and let $(\sigma, W)$ be a subrepresentation of $\rho_{\lambda}\otimes\overline{\mathbb{Q}}_{\ell}$. Then $\underline{G}_{\ell}$ preserves $\overline{W}^{\mathrm{ss}}$. Denote by $\underline{G}_W$ the image of $\underline{G}_{\ell}$ in $\mathrm{GL}_{\overline{W}^{\mathrm{ss}}}$. We call $\underline{G}_W$ the algebraic envelope of $W$.
	\end{defn/thm}
	
	\begin{prop}\label{Manganese}
		Under the hypotheses of the above theorem, for all but finitely many $\lambda$ and each subrepresentation $(\sigma,W)$ of $\rho_{\lambda}\otimes\overline{\mathbb{Q}}_{\ell}$, we have
		$$\dim\underline{G}_W\le\dim\mathbf{G}_W^{\circ},$$
		where $\mathbf{G}_W^{\circ}$ is the connected monodromy group of $\sigma$.
	\end{prop}
	\begin{proof}
		Put $r=\dim W$. By \cite[Proposition 3.9(i)-(iv)]{Hu23b}, there exist connected algebraic groups
		$$\underline{G}_W'\subseteq\underline{H}_W'\subseteq\mathrm{GL}_{r,\overline{\mathbb{F}}_{\ell}}$$
		such that $\underline{G}_W$ is isomorphic to the reductive quotient of $\underline{G}_W'$, and $\mathbf{G}_W^{\circ}$ and $\underline{H}_W'$ are isomorphic to the generic and special fibers, respectively, of a common split reductive group scheme. Hence
		$$\dim\underline{G}_W\le\dim\underline{G}_W'\le\dim\underline{H}_W'=\dim\mathbf{G}_W^{\circ}.$$
		This gives the conclusion.
	\end{proof}
	
	\begin{thm}\cite[Theorem 4.1]{Hu23b}\label{Iron}
		Let $K$ be imaginary CM or totally real and let $\pi$ be a regular algebraic polarized cuspidal automorphic representation of $\mathrm{GL}_n(\mathbb{A}_K)$. Then the associated compatible system $\{\rho_{\pi,\lambda}\}_{\lambda}$ as in \autoref{Titanium} satisfies conditions (a) and (b) in \autoref{Vanadium}.
	\end{thm}
	
	\subsection{Polarization}
	
	\begin{defn}
		Let $K$ be a CM or totally real field and let $K^+$ be its maximal totally real subfield. Let $\rho:\mathrm{Gal}_K\to\mathrm{GL}(V)$ be an $\ell$-adic representation over $\overline{\mathbb{Q}}_{\ell}$ and let $\mu:\mathrm{Gal}_{K^+}\to\overline{\mathbb{Q}}_{\ell}^{\times}$ be a continuous character. The pair $(\rho,\mu)$ is called polarized if for every infinite place $v$ of $K^+$, there exist $\epsilon_v\in\{\pm1\}$ and a nondegenerate bilinear form $\langle\cdot,\cdot\rangle_v$ on $V$ satisfying
		\begin{align*}
			\langle x,y\rangle_v&=\epsilon_v\langle y,x\rangle_v,\\
			\langle\rho(g)x,\rho(c_vgc_v)y\rangle_v&=\mu(g)\langle x,y\rangle_v
		\end{align*}
		for all $x,y\in V$ and $g\in\mathrm{Gal}_K$, where $c_v\in\mathrm{Gal}_{K^+}$ is a complex conjugation at $v$. If $K$ is imaginary CM, we additionally require $\epsilon_v=-\mu(c_v)$. A polarized pair $(\rho,\mu)$ is called totally odd if $\epsilon_v=1$ for every infinite place $v$ of $K^+$.
	\end{defn}
	
	The following result \autoref{Carbon} is a slightly stronger variant of \cite[Proposition 5.10]{FW25}, whose proof relies on \cite[Proposition 17]{Xi19}. For the reader's convenience, we include a proof, using a different auxiliary representation and the results of Larsen and Pink \cite[proof of Proposition 4.11 and Lemma 6.11]{LP92}.
	
	\begin{lem}\label{Sodium}
		Let $\rho$ be a semisimple $\lambda$-adic Galois representation of a number field $K$ and let $\mathbf{G}$ be its algebraic monodromy group.
		\begin{enumerate}
			\item[(i)] The connected monodromy group of $\mathrm{End}(\rho)$ is $\mathbf{G}^{\circ}/\left(\mathbf{G}^{\circ}\cap\mathbb{G}_m\right)$.
			\item[(ii)] If $K$ is imaginary CM and $\rho$ is polarized. Let $K^+$ be its maximal totally real subfield. Then the derived monodromy group of $\mathrm{Ind}_K^{K^+}\rho$ is isomorphic to $\mathbf{G}^{\mathrm{der}}$.
			\item[(iii)] If $K$ is imaginary CM and $\rho$ is polarized. Let $K^+$ be its maximal totally real subfield. Then the connected monodromy group of $\mathrm{Ind}_K^{K^+}\mathrm{End}(\rho)$ is isomorphic to $\mathbf{G}^{\circ}/(\mathbf{G}^{\circ}\cap\mathbb{G}_m)$.
		\end{enumerate}
	\end{lem}
	\begin{proof}
		For (i), the conjugation representation $\mathrm{Ad}:\mathrm{GL}_n\to\mathrm{GL}_{n^2}$ has kernel $\mathbb{G}_mI_n$. Since $\mathrm{End}(\rho)=\mathrm{Ad}\circ\rho$ and $\rho(\mathrm{Gal}_K)$ is Zariski dense in $\mathbf{G}$, the algebraic monodromy group of $\mathrm{End}(\rho)$ is $\mathrm{Ad}(\mathbf{G})$. Taking identity components yields that the connected monodromy group of $\mathrm{End}(\rho)$ is
		$$\mathrm{Ad}(\mathbf{G}^{\circ})\simeq\mathbf{G}^{\circ}/\left(\mathbf{G}^{\circ}\cap\mathbb{G}_m\right).$$
		For (ii), put $\varphi=\mathrm{Ind}_K^{K^+}\rho$. By polarization of $\rho$,
		$$\varphi|_{\mathrm{Gal}_K}\simeq\rho\oplus\left(\rho^{\vee}\otimes\mu|_{\mathrm{Gal}_K}\right).$$
		Put $\psi=\rho\oplus\mu|_{\mathrm{Gal}_K}$. Define a homomorphism
		$$\iota:\mathrm{GL}_n\times\mathbb{G}_m\to\mathrm{GL}_{2n},\quad(x,z)\mapsto\diag\{x,zx^{-t}\}.$$
		Then $\varphi|_{\mathrm{Gal}_K}=\iota\circ\psi$. Let $\mathbf{H}$ be the algebraic monodromy group of $\psi$, viewed as a subgroup of $\mathrm{GL}_n\times\mathbb{G}_m$. Projection onto the first factor is a surjection from $\mathbf{H}^{\circ}$ to $\mathbf{G}^{\circ}$. In particular, it maps $\mathbf{H}^{\mathrm{der}}$ onto $\mathbf{G}^{\mathrm{der}}$. Since every character is trivial on a derived group,
		\begin{equation}\label{3}
			\mathbf{H}^{\mathrm{der}}=\mathbf{G}^{\mathrm{der}}\times1\simeq\mathbf{G}^{\mathrm{der}}.
		\end{equation}
		Since $\mathrm{Gal}_K$ has finite index in $\mathrm{Gal}_{K^+}$, restriction to $\mathrm{Gal}_K$ does not change the connected monodromy group of $\varphi$. Hence the derived monodromy of $\varphi$ is isomorphic to $\iota(\mathbf{H}^{\mathrm{der}})$, which is isomorphic to $\mathbf{G}^{\mathrm{der}}$ by (\ref{3}).
		
		For (iii), by polarization, $\mathrm{End}(\rho)^c\simeq\mathrm{End}(\rho)$. Thus the restriction of $\mathrm{Ind}_K^{K^+}\mathrm{End}(\rho)$ to $\mathrm{Gal}_K$ is $\mathrm{End}(\rho)^{\oplus 2}$. Since restriction to a finite-index subgroup does not change the connected monodromy group, the conclusion follows from (i).
	\end{proof}
	
	\begin{prop}\label{Carbon}
		Let $K$ be imaginary CM or totally real and let $K^+$ be its maximal totally real subfield. Let 
		$$\{\rho_{\lambda}:\mathrm{Gal}_K\to\mathrm{GL}_n(\overline{E}_{\lambda})\}_{\lambda}$$
		be an $E$-rational polarized semisimple Serre compatible system of $K$. Denote by $\mathbf{G}_{\lambda}$ the algebraic monodromy group of $\rho_{\lambda}$. Fix a complex conjugation $c\in\mathrm{Gal}_{K^+}$, and put $c_{\lambda}=\mathrm{Ind}_K^{K^+}(\rho_{\lambda})(c)$, which acts by conjugation on $\mathbf{G}_{\lambda}^{\mathrm{der}}$ via the identification of the derived monodromy groups of $\rho_{\lambda}$ and $\mathrm{Ind}_K^{K^+}(\rho_{\lambda})$ by \autoref{Sodium}(ii). Assume that at some coefficient prime $\lambda_0$, the representation $\rho_{\lambda_0}$ is Lie-irreducible and that
		\begin{equation}
			c_{\lambda_0}\text{ acts on }\mathbf{G}_{\lambda_0}^{\mathrm{der}}\text{ by an inner automorphism.}
			\tag{*}
		\end{equation}
		Then:
		\begin{enumerate}
			\item[(i)] $c_{\lambda}$ acts on $\mathbf{G}_{\lambda}^{\mathrm{der}}$ by an inner automorphism for every $\lambda$.
			\item[(ii)] If $\rho_{\lambda_0}$ has multiplicity-free derived formal character, then for every coefficient prime $\lambda$, any irreducible constituent $U\subseteq\rho_{\lambda}$ such that $U^{\mathrm{der}}$ contains a self-dual irreducible summand is polarized. Moreover, if $\rho_{\lambda}$ is totally odd, then $U$ is also totally odd.
			\item[(iii)] If $K$ is totally real, then (ii) holds without condition $(\ast)$.
		\end{enumerate}
	\end{prop}
	\begin{proof}[Proof of \autoref{Carbon}(i)]
		Put $\varphi_{\lambda}=\mathrm{Ind}_K^{K^+}\mathrm{End}(\rho_{\lambda})$ and let $\mathbf{H}_{\lambda}$ be its algebraic monodromy group. The family $\{\varphi_{\lambda}\}_{\lambda}$ is a semisimple $E$-rational compatible system. By \autoref{Sodium}(iii), $$\mathbf{H}_{\lambda}^{\circ}\simeq\mathbf{G}_{\lambda}^{\circ}/(\mathbf{G}_{\lambda}^{\circ}\cap\mathbb{G}_m).$$
		Since $\mathbf{G}_{\lambda_0}^{\circ}$ is connected reductive, we have
		$$\mathbf{G}_{\lambda_0}^{\circ}=Z(\mathbf{G}_{\lambda_0}^{\circ})^{\circ}\mathbf{G}_{\lambda_0}^{\mathrm{der}}.$$
		Since $\rho_{\lambda_0}$ is Lie-irreducible, Schur's lemma implies that $Z(\mathbf{G}_{\lambda_0}^{\circ})=\mathbf{G}_{\lambda_0}^{\circ}\cap\mathbb{G}_m$. Consequently,
		$$\mathbf{H}_{\lambda_0}^{\circ}\simeq\mathbf{G}_{\lambda_0}^{\circ}/(\mathbf{G}_{\lambda_0}^{\circ}\cap\mathbb{G}_m)\simeq\mathbf{G}_{\lambda_0}^{\mathrm{der}}/(\mathbf{G}_{\lambda_0}^{\mathrm{der}}\cap\mathbb{G}_m)$$
		As a central quotient of $\mathbf{G}_{\lambda_0}^{\mathrm{der}}$, $\mathbf{H}_{\lambda_0}^{\circ}$ is semisimple. Hence \autoref{Hydrogen}(ii) implies that each $\mathrm{H}_{\lambda}^{\circ}$ is semisimple.
		
		Write $C_{\lambda}=\varphi_{\lambda}(c)$ and $D_{\lambda}=C_{\lambda}H_{\lambda}^{\circ}$. Under the natural central quotient $q_{\lambda}:\mathbf{G}_{\lambda}^{\mathrm{der}}\twoheadrightarrow\mathbf{H}_{\lambda}^{\mathrm{der}}=\mathbf{H}_{\lambda}^{\circ}$, the conjugation actions of $c_{\lambda}$ and $C_{\lambda}$ are compatible:
		\begin{equation}\label{4}
			q_{\lambda}(c_{\lambda}xc_{\lambda}^{-1})=C_{\lambda}q_{\lambda}(x)C_{\lambda}^{-1},\quad\forall x\in\mathbf{G}_{\lambda}^{\mathrm{der}}.
		\end{equation}
		Let $L/K^+$ be the connectedness field of $\{\varphi_{\lambda}\}_{\lambda}$ provided by \autoref{Hydrogen}(i). Then $D_{\lambda}$ is the Zariski closure of $\varphi_{\lambda}(c\mathrm{Gal}_L)$. Denote by $\mathrm{ch}$ the characteristic polynomial map $\mathrm{GL}_n\to\mathbb{A}^{n-1}\times\mathbb{G}_m$, with $n=\dim\varphi_{\lambda}$. By \cite[Proposition 4.9 and Lemma 6.11]{LP92}, each $\mathrm{ch}(D_{\lambda})$ is Zariski closed and the dimension of $\mathrm{ch}(D_{\lambda})$ is independent of $\lambda$.
		
		By condition $(\ast)$ and (\ref{4}), there exists $h\in\mathbf{H}_{\lambda_0}^{\circ}$ such that 
		$$C_{\lambda_0}xC_{\lambda_0}^{-1}=hxh^{-1},\quad\forall x\in\mathbf{H}_{\lambda_0}^{\circ}.$$
		Put $z=h^{-1}C_{\lambda_0}$. Then $z$ centralizes $\mathbf{H}_{\lambda_0}^{\circ}$ and
		$$D_{\lambda_0}=C_{\lambda_0}\mathbf{H}_{\lambda_0}^{\circ}=z\mathbf{H}_{\lambda_0}^{\circ}.$$
		Choose a maximal torus $T\subseteq\mathbf{H}_{\lambda_0}^{\circ}$. Since $C_{\lambda_0}^2=1$ and $C_{\lambda_0}$ commutes with $h$, we have $z^2=h^{-2}\in Z(\mathbf{H}_{\lambda_0}^{\circ})$. This center is finite since $\mathbf{H}_{\lambda_0}^{\circ}$ is semisimple, so $z$ has finite order. Put $H=\langle z,T\rangle$. Since $z$ centralizes $T$ and has finite order, $H$ is an algebraic group with $H^{\circ}=T$. By \cite[Proposition 4.9]{LP92}, $\mathrm{ch}(zT)$ is Zariski closed. Since the conjugates of $T$ are Zariski dense in $\mathbf{H}_{\lambda_0}^{\circ}$ and $z$ centralizes $\mathbf{H}_{\lambda_0}^{\circ}$, the conjugates of $zT$ are Zariski dense in $D_{\lambda_0}$. Since characteristic polynomials are invariant under conjugation, we have
		\begin{equation}\label{5}
			\mathrm{ch}(D_{\lambda_0})=\mathrm{ch}(zT)
		\end{equation}
		Moreover, $z$ is semisimple and centralizes $T$, so all elements in $zT$ are simultaneously diagonalizable. Consequently, the restriction of $\mathrm{ch}$ to $zT$ has finite fibers. Hence
		$$\dim(\mathrm{ch}(zT))=\dim(zT)=\dim T=\rk\mathbf{H}_{\lambda_0}^{\circ}$$
		Combining this with (\ref{5}), the $\lambda$-independence of $\dim\mathrm{ch}(D_{\lambda})$ and the $\lambda$-independence of $\rk\mathbf{H}_{\lambda}^{\circ}$, we obtain
		$$\dim\mathrm{ch}(D_{\lambda})=\dim\mathrm{ch}(D_{\lambda_0})=\rk\mathbf{H}_{\lambda_0}^{\circ}=\rk\mathbf{H}_{\lambda}^{\circ}.$$
		The implication proved in \cite[proof of Proposition 4.11]{LP92} therefore implies that $C_{\lambda}$ acts on $\mathbf{H}_{\lambda}^{\circ}$ by an inner automorphism for every $\lambda$. Let $h=h_{\lambda}\in\mathbf{H}_{\lambda}^{\circ}$ be an element satisfying
		$$C_{\lambda}xC_{\lambda}^{-1}=hxh^{-1},\quad\forall x\in\mathbf{H}_{\lambda}^{\circ}.$$ 
		Since $q_{\lambda}$ is surjective, choose $g\in\mathbf{G}_{\lambda}^{\mathrm{der}}$ with $q_{\lambda}(g)=h$. Then (\ref{4}) implies that
		$$q_{\lambda}(c_{\lambda}yc_{\lambda}^{-1})=q_{\lambda}(gyg^{-1}),\quad\forall y\in\mathbf{G}_{\lambda}^{\mathrm{der}}.$$
		Hence we have a morphism
		$$\mathbf{G}_{\lambda}^{\mathrm{der}}\to\ker q_{\lambda}\qquad y\mapsto gy^{-1}g^{-1}c_{\lambda}yc_{\lambda}^{-1}.$$
		Since $\ker q_{\lambda}$ is finite and $\mathbf{G}_{\lambda}^{\mathrm{der}}$ is connected, this morphism is constant. Evaluating at $y=1$ shows that it is identically equal to 1. Thus
		$$c_{\lambda}yc_{\lambda}^{-1}=gyg^{-1},\quad\forall y\in\mathbf{G}_{\lambda}^{\mathrm{der}}.$$
		Hence conjugation by $c_{\lambda}$ on $\mathbf{G}_{\lambda}^{\mathrm{der}}$ is inner. This completes the proof of (i).
	\end{proof}
	
	\begin{proof}[Proof of \autoref{Carbon}(ii) and (iii)]
		Let $w$ be a weight in the self-dual irreducible summand of $U^{\mathrm{der}}$. Then both $w$ and $-w$ are weights of $U^{\mathrm{der}}$. By \autoref{Carbon}(i), $-w$ is a weight of $(U^c)^{\mathrm{der}}$. On the other hand, since twist by a character does not change the derived formal character, $-w$ is a weight of $\left(U^{\vee}\otimes\mu_{\lambda}|_{\mathrm{Gal}_K}\right)^{\mathrm{der}}$. By polarization,
		$$\rho_{\lambda}^c\simeq\rho_{\lambda}^{\vee}\otimes\mu_{\lambda}|_{\mathrm{Gal}_K}.$$
		By \autoref{Hydrogen}(ii), the derived formal character of $\rho_{\lambda}$, and hence that of its dual twist, is multiplicity-free. Since $U^c$ and $U^{\vee}\otimes\mu_{\lambda}|_{\mathrm{Gal}_K}$ are irreducible constituents of this common representation and share the weight $-w$ in the derived formal character, we conclude that
		$$U^c\simeq U^{\vee}\otimes\mu_{\lambda}|_{\mathrm{Gal}_K}.$$
		
		For every complex conjugation $c_v$ at an infinite place $v$ of $K^+$, we have $c_vc^{-1}\in\mathrm{Gal}_K$. Thus $U$ is polarized. Moreover, $U$ occurs only once in $\rho_{\lambda}$, since the derived formal character of $\rho_{\lambda}$ is multiplicity-free. Hence, by Schur's lemma, the nondegenerate polarization pairing on $\rho_{\lambda}$ restricts nondegenerately to $U$. The restricted pairing has the same symmetric or alternating type and multiplier as the original pairing. Consequently, if $\rho_{\lambda}$ is totally odd, then so is $U$. This proves (ii).
		
		If $K$ is totally real, since $c\in\mathrm{Gal}_K$, we have $U^c\simeq U$, and hence their derived formal character coincide without (i). Therefore the proof of (ii) above applies without condition $(\ast)$. This proves (iii).
	\end{proof}
	
	\subsection{Irreducibility of automorphic Galois representations}
	
	Let $K$ be a totally real or imaginary CM field and let $\pi$ be a regular algebraic polarized cuspidal automorphic representation of $\mathrm{GL}_n(\mathbb{A}_K)$. Denote by $E_\pi$ the fixed field of
	$$\{\gamma\in\mathrm{Aut}(\mathbb{C}):{}^\gamma\!\pi^\infty\simeq\pi^\infty\}.$$
	By \autoref{Titanium}(ii), \cite[Theorem 3.13]{Cl90}, and strong multiplicity one, $E_\pi$ is precisely the minimal coefficient field of the compatible system associated to $\pi$, namely the number field generated over $\mathbb{Q}$ by the coefficients of the Frobenius characteristic polynomials $P_v(T)$ at almost all finite places $v$ of $K$. Moreover, by \cite[Lemma 1.1]{PT15}, the coefficient field $E$ in \autoref{Titanium} may be taken to be $E_\pi$.
	\begin{thm}\cite[Theorem 1.7]{PT15}\label{Boron}
		There exist a finite CM extension $E/E_\pi$ and a set $\mathcal{L}$ of rational primes of Dirichlet density one such that for every $\ell\in\mathcal{L}$ and every conjugation-invariant place $\lambda$ of $E$ above $\ell$, the representation $\rho_{\pi,\lambda|_{E_\pi}}\otimes\overline{\mathbb{Q}}_\ell$ is irreducible.
		In particular, there exists a set $\mathcal{L}'$ of rational primes of positive Dirichlet density such that $\rho_{\pi,\lambda}\otimes\overline{\mathbb{Q}}_\ell$ is irreducible for every place $\lambda$ of $E_\pi$ above $\ell\in\mathcal{L}'$.
	\end{thm}
	
	\begin{prop}\label{Magnesium}
		Let $K$ be totally real or imaginary CM, and let $\{\rho_{\lambda}\}_{\lambda}$ be the $E$-rational compatible system attached to a regular algebraic polarized cuspidal automorphic representation $(\pi,\chi)$ of $\mathrm{GL}_n(\mathbb{A}_K)$. Then the set of rational primes $\ell$ for which there exists a place $\lambda|\ell$ such that $\rho_{\lambda}$ is reducible and every irreducible constituent of $\rho_{\lambda}$ is polarized and totally odd has Dirichlet density zero.
	\end{prop}
	\begin{proof}
		By \cite[Proposition 5.3.2]{BLGGT14}, for a set $\mathcal{L}$ of rational primes of Dirichlet density one, every irreducible constituent $\sigma$ of $\rho_\lambda$ with $\lambda\mid\ell\in\mathcal{L}$ satisfies that $\overline{\sigma}^{\mathrm{ss}}|_{\mathrm{Gal}_{K(\zeta_\ell)}}$ is irreducible. After removing finitely many primes from $\mathcal{L}$, if every irreducible constituent of $\rho_\lambda$ is polarized and totally odd, then \cite[Theorem 4.10]{Hu23b} implies that these constituents become irreducible and automorphic over a common finite Galois extension of $K$. By \cite[Proposition 4.9]{Hu23b}, $\rho_\lambda$ is irreducible, as desired.
	\end{proof}
	
	\begin{cor}\label{Helium}
		Let $K$ be totally real or imaginary CM, and let $\{\rho_{\lambda}\}_{\lambda}$ be the $E$-rational compatible system attached to a regular algebraic polarized cuspidal automorphic representation $(\pi,\chi)$ of $\mathrm{GL}_n(\mathbb{A}_K)$. Then the set of rational primes $\ell$ for which there exists a place $\lambda|\ell$ such that $\rho_{\lambda}$ is reducible and the irreducible constituents of $\rho_{\lambda}$ have pairwise distinct dimensions has Dirichlet density zero.
	\end{cor}
	\begin{proof}
		Suppose that the irreducible constituents of $\rho_\lambda$ have pairwise distinct dimensions. By Schur's lemma, each polarization pairing of $\rho_\lambda$ vanishes on any two distinct irreducible constituents. Therefore, its restriction to each irreducible constituent is nondegenerate, with the same symmetry and multiplier. Thus every irreducible constituent is polarized and totally odd, and the assertion follows from \autoref{Magnesium}.
	\end{proof}
	
	\section{Proof of the Main Theorem}
	
	\begin{lem}\label{Nitrogen}
		If the tautological representation of $\mathbf{G}_{\lambda_0}^{\mathrm{der}}$ is the standard representation of the exceptional group $G_2$, then for each $\lambda$, the tautological representation of $\mathbf{G}_{\lambda}^{\mathrm{der}}$ is one of the following.
		\begin{enumerate}
			\item[(i)] $\left(G_2,\mathrm{Std}\right)$.
			\item[(ii)] $\left(\mathrm{SL}_3,\mathrm{Std}\oplus\mathrm{Std}^{\vee}\oplus\mathbf{1}\right)$.
			\item[(iii)] $\left(\mathrm{SO}_4,\mathrm{Std}\oplus(\mathrm{Sym}^2(S^+))\right)$, where $S^+$ is one of the half-spin representations of $\mathrm{Spin}_4$, and its symmetric square is viewed as a representation of $\mathrm{SO}_4$.
		\end{enumerate}
	\end{lem}
	\begin{proof}
		The formal character of $\mathbf{G}_{\lambda_0}^{\mathrm{der}}$ is
		$$\{x,-x,y,-y,x+y,-x-y,0\},$$
		where $x$ and $y$ are linearly independent characters. By \autoref{Hydrogen}, all $\mathbf{G}_{\lambda}^{\mathrm{der}}$ have the same formal character as above. Equip the weight space with the inner product for which $(x,x)=(y,y)=1$ and $(x,y)=-1/2$. The six nonzero weights form the vertices of a regular hexagon, with the zero weight at its center. Since every root is a difference of two weights, the possible roots are
		$$\begin{array}{c|l}
			\text{Length} & \text{Candidate roots} \\ \hline
			1
			& \pm x,\ \pm y,\ \pm(x+y) \\[2pt]
			\sqrt{3}
			& \pm(x-y),\ \pm(2x+y),\ \pm(x+2y) \\[2pt]
			2
			& \pm 2x,\ \pm 2y,\ \pm 2(x+y)
		\end{array}$$
		
		If the length-two candidates were roots, the corresponding coroot pairings would give
		$$\langle y,(2x)^{\vee}\rangle=-\frac{1}{2},\quad\langle x,(2y)^{\vee}\rangle=-\frac{1}{2},\quad\langle x,(2(x+y))^{\vee}\rangle=\frac{1}{2}.$$
		None is an integer. Thus the length-two candidates cannot be roots. The remaining twelve candidates are exactly the six short roots and six long roots of $G_2$.
		
		If only one length occurs, rank two gives two nonparallel roots, whose reflections generate all six roots of that length. If a nonperpendicular pair of short root and long root occurs, its reflections generate a $60^{\circ}$ rotation, giving all twelve possible roots. Otherwise, if both lengths occur, and every short root is perpendicular to every long root, then there is exactly one pair of opposite roots for each length. We distinguish three cases.
		
		\begin{itemize}
			\item All twelve possible roots occur in $\mathbf{G}_{\lambda}^{\mathrm{der}}$. This is case (i).
			\item Exactly the six roots of one length occur. If they are the short roots, they form a root system of type $A_2$. Choose simple roots $\alpha$ and $\beta$. Then $\theta=\alpha+\beta$ is a highest weight in the formal character of $\mathbf{G}_{\lambda}^{\mathrm{der}}$, which gives an irreducible constituent of the seven-dimensional representation. But this constituent is the eight-dimensional adjoint representation of $A_2$, a contradiction. Hence the six long roots occur. Choose $x-y$ and $x+2y$ as simple roots. Then there are three highest weights $x,x+y,0$ among the formal character of $\mathbf{G}_{\lambda}^{\mathrm{der}}$. The highest weights $x$ and $x+y$ correspond the standard representation of $\mathrm{SL}_3$ and its dual, respectively. The remaining weight zero gives the trivial summand. This gives case (ii).
			\item A short root $\alpha$ and a long root $\beta$ occur, but not all twelve roots occur. Then $\alpha$ and $\beta$ are perpendicular. The set of weights is
			$$\{0,\pm\alpha,\frac{\pm\alpha\pm\beta}{2}\}.$$ 
			Choose $\alpha,\beta$ as simple roots. Then the highest weights are $\alpha,\frac{\alpha+\beta}{2}$. On the simply connected cover $\mathrm{SL}_2\times\mathrm{SL}_2$, the corresponding summands are $\mathrm{Sym}^2(\mathrm{Std})\otimes\mathbf{1}$ and $\mathrm{Std}\otimes\mathrm{Std}$. The kernel of the latter one is $\mu_2=\{(I,I),(-I,-I)\}$, which also acts trivially on the former. Thus their direct sum descends to a faithful representation of $\mathrm{SO}_4$. This gives case (iii).
		\end{itemize}
	\end{proof}
	
	Since $\mathbf{G}_{\lambda}$ normalizes $\mathbf{G}_{\lambda}^{\mathrm{der}}$, it permutes the irreducible summands of the derived representation. In case (iii), the two irreducible constituents of $\rho_{\lambda}$ have dimensions 4 and 3, which is different, so each is $\mathbf{G}_{\lambda}$-stable and therefore defines an irreducible constituent of $\rho_{\lambda}$. Hence by \autoref{Helium}, the set of rational primes $\ell$ for which case (iii) occurs for at least one $\lambda|\ell$ has Dirichlet density zero. To deal with case (ii), we first make the following normalization.
	\begin{lem}[Normalization of automorphic compatible system of type $G_2$]\label{Chlorine}
		Let $K$ be imaginary CM or totally real, and let $\pi$ be a regular algebraic polarized cuspidal automorphic representation of $\mathrm{GL}_n(\mathbb{A}_K)$ with associated $E$-rational compatible system $\{\rho_{\lambda}\}_{\lambda}$. Suppose that there is a coefficient prime $\lambda_0$ such that $\mathbf{G}_{\lambda_0}^{\mathrm{der}}$ acts via the standard seven-dimensional irreducible representation of $G_2$. Then, after twisting $\pi$ by an algebraic Hecke character of $K$ and enlarging $E$ if necessary, we may assume that $\mathbf{G}_{\lambda_0}$ is of type $G_2$, and that $\mathbf{G}_{\lambda}$ is connected and semisimple for every coefficient prime $\lambda$.
	\end{lem}
	\begin{proof}
		Since $\mathbf{G}_{\lambda_0}$ normalizes $\mathbf{G}_{\lambda_0}^{\mathrm{der}}=G_2$, every automorphism of $G_2$ is inner, and the centralizer of $G_2$ in $\mathrm{GL}_7$ consists of scalars, we obtain
		$$G_2\subseteq\mathbf{G}_{\lambda_0}\subseteq\mathrm{N}_{\mathrm{GL}_7}(G_2)=\mathbb{G}_mG_2.$$
		Since $\mathbb{G}_m$ and $G_2$ intersect trivially, we may write $\rho_{\lambda_0}=\alpha r$, with $\alpha$ a character and $r$ valued in $G_2$. Let $\beta$ be the unique one-dimensional subrepresentation of $\wedge^3(\rho_{\lambda_0})$. Then $\alpha=\det(\rho_{\lambda_0})\beta^{-2}$, so $\alpha$ is de Rham and it corresponds to an algebraic Hecke character $\chi$. Denote by $\{\alpha_{\lambda}\}_{\lambda}$ the compatible system attached to $\chi$ after possibly enlarging the coefficient field. Replace $\pi$ with $\pi\otimes(\chi^{-1}\circ\det)$, we may assume that the associated compatible system is $\{\rho_{\lambda}\alpha_{\lambda}^{-1}\}_{\lambda}$. In particular, $\mathbf{G}_{\lambda_0}=G_2$. Since this group is connected, \autoref{Hydrogen}(i) implies that each $\mathbf{G}_{\lambda}$ is connected. By \autoref{Hydrogen}(ii), both the rank and the semisimple rank of $\mathbf{G}_{\lambda}$ are independent of $\lambda$. These ranks agree at $\lambda_0$, so they agree at every $\lambda$. Thus each $\mathbf{G}_{\lambda}$ is semisimple. 
	\end{proof}
	
	Twisting by a character preserves irreducibility, so it suffices to work with this new $\pi$ and its associated system. From now on, we write $\pi$ for the twisted automorphic representation, $\{\rho_{\lambda}\}_{\lambda}$ for its associated compatible system, and $\mathbf{G}_{\lambda}$ for the algebraic monodromy group of $\rho_{\lambda}$.
	
	Since the standard representation of $G_2$ has a zero weight, by \autoref{Hydrogen}(ii), the formal character of each $\mathbf{G}_{\lambda}$ also has a zero weight. Since $\mathbf{G}_{\lambda}$ is connected, for any $g\in\mathbf{G}_{\lambda}(E_{\lambda})$, its semisimple part lies in a maximal torus over $\overline{E}_{\lambda}$ and has the same eigenvalues as $g$. Hence $1$ is an eigenvalue of $g$. Write its characteristic polynomial as
	$$P_g(T)=(T-1)\left(T^6-AT^5+BT^4+\cdots\right)\in E_{\lambda}[T],$$
	and define
	$$\Delta(g):=A^2+4A-4B\in E_{\lambda}.$$
	Denote by $S$ the exceptional set of $\{\rho_{\lambda}\}_{\lambda}$. For each $v\in\Sigma_K\backslash S$, define
	$$\Delta_v:=\Delta(\rho_{\lambda}(\mathrm{Frob}_v))\in E,$$
	where $\lambda$ and $v$ have different residue characteristics. This is independent of the choice of such $\lambda$.
	
	Denote by $\mathcal{L}$ the set of places $\lambda$ in case (ii). Since $\mathbf{G}_{\lambda}$ is connected and semisimple, we have
	$$\rho_{\lambda}=U_{\lambda}\oplus U_{\lambda}^{\vee}\oplus\mathbf{1},\quad\forall\lambda\in\mathcal{L}.$$
	For each $\lambda\in\mathcal{L}$ and each $v\in\Sigma_K\backslash S$ such that they have different residue characteristics, define
	$$a_{\lambda,v}=\mathrm{tr}\left(U_{\lambda}(\mathrm{Frob}_v)\right),\quad b_{\lambda,v}=\mathrm{tr}\left(U_{\lambda}^{\vee}(\mathrm{Frob}_v)\right)\in\overline{E}_{\lambda}.$$
    The characteristic polynomials of $U_{\lambda}(\mathrm{Frob}_v)$ and $U_{\lambda}^{\vee}(\mathrm{Frob}_v)$ are
	$$P_{U_{\lambda}(\mathrm{Frob}_v)}=T^3-a_{\lambda,v}T^2+b_{\lambda,v}T-1,\quad P_{U_{\lambda}^{\vee}(\mathrm{Frob}_v)}=T^3-b_{\lambda,v}T^2+a_{\lambda,v}T-1.$$
	Hence we have
	\begin{equation}\label{1}
	\Delta_v=(a_{\lambda,v}-b_{\lambda,v})^2\in E
	\end{equation}
	
	For each $v\in\Sigma_K\backslash S$, denote by $P_v(T)\in E[T]$ the common characteristic polynomial of $\rho_{\lambda}(\mathrm{Frob}_v)$ such that $\lambda$ and $v$ have different residual characteristics. Let $E_{\pi}$ be the fixed field of the subgroup of $\mathrm{Aut}(\mathbb{C})$ consisting of elements $\sigma$ such that $^{\sigma}\pi^{\infty}\simeq\pi^{\infty}$. As in \cite[p. 10]{PT15}, we take $E=E_{\pi}$ as the coefficient field of the compatible system associated to $\pi$. By strong multiplicity one, $E$ is the number field generated by the coefficients of $P_v(T)$ for $v\not\in S$.
	\begin{lem}\label{Beryllium}
		Let $M$ be the Galois closure of $E/\mathbb{Q}$. Put $X=\mathrm{Hom}_{\mathbb{Q}}(E,M)$. 
		\begin{enumerate}
			\item[(i)] There exist infinitely many odd rational primes $p$ such that $p$ splits completely in $M$ and that $\mathbf{G}_{\lambda}$ is of type $G_2$ for every $\lambda|p$.
			\item[(ii)] Fix a rational prime $p$ as in (i). Then for each $\lambda|p$, up to a change of coordinates, $\rho_{\lambda}$ factors through $G_2(\mathbb{Q}_p)$. With these choices of coordinates, the product representation
			$$\rho_p:=(\rho_{\lambda})_{\lambda|p}:\mathrm{Gal}_K\to\prod_{\lambda|p}G_2(\mathbb{Q}_p)$$
			has open image.
			\item[(iii)] For each integer $N\ge1$, there exist primes $v_0,v_1,\cdots,v_N\in\Sigma_K\backslash S$, such that the $(N+1)[E:\mathbb{Q}]$ classes
			$$[\sigma(\Delta_{v_i})],\quad0\le i\le N,\sigma\in X$$
			are linearly independent over $\mathbb{F}_2$ in $M^{\times}/(M^{\times})^2$. In particular, for $r_i=\Delta_{v_i}/\Delta_{v_0}$, the classes
			$$[\sigma(r_i)],\quad1\le i\le N,\sigma\in X$$
			are linearly independent over $\mathbb{F}_2$ in $M^{\times}/(M^{\times})^2$.
		\end{enumerate}
	\end{lem}
	
	We use the following consequence of \cite{Pi98} to prove openness in part (ii).
	\begin{lem}\label{Oxygen}
		Let $F$ be a $p$-adic field, and let $\Gamma\subseteq G_2(F)$ be a compact Zariski-dense subgroup. If $\mathrm{tr}(\mathrm{Ad}(\gamma))\in\mathbb{Q}_p$ for every $\gamma\in\Gamma$, then up to conjugacy, $\Gamma$ is an open subgroup of $G_2(\mathbb{Q}_p)$.
	\end{lem}
	\begin{proof}
		Since $G_2/F$ is an absolutely simple connected adjoint group with irreducible adjoint representation, \cite[Theorem 0.7]{Pi98} gives a model $H$ of $G_2$ over a closed subfield $E\subseteq F$, such that $\Gamma$ is open in $H(E)$. By \cite[Proposition 0.6(a)]{Pi98}, $E$ is the fraction field of the closed subring of $F$ generated by all the traces $\mathrm{tr}(\mathrm{Ad}(\gamma))$ for $\gamma\in\Gamma$. This subring lies in $\mathbb{Q}_p$ by assumption. Hence the fraction field $E=\mathbb{Q}_p$. Every $\mathbb{Q}_p$-form of $G_2/F$ is split, so $H\simeq G_2$ over $\mathbb{Q}_p$. Thus, up to conjugation, $\Gamma$ is an open subgroup of $G_2(\mathbb{Q}_p)$.
	\end{proof}
	
	\begin{proof}[Proof of \autoref{Beryllium}(i)]
		We first show that $E$ and $M$ are totally real. For each $v\in\Sigma_K\backslash S$, purity implies that, as $\kappa$ runs over the embeddings $E\hookrightarrow\mathbb{C}$, all roots of all the polynomials $\kappa(P_v)(T)$ have the same absolute value. Since $\mathbf{G}_{\lambda}$ is semisimple in our situation, $\det\rho_{\lambda}$ is trivial. Hence this common absolute value is 1. Moreover, for an irreducible $\rho_{\lambda_0}$, the algebraic monodromy group acts as the standard seven-dimensional representation of $G_2$, which is self-dual. Hence $\rho_{\lambda_0}\simeq\rho_{\lambda_0}^{\vee}$. By compatibility and Chebotarev, the roots of each $P_v(T)$ are stable under inversion, with multiplicities. Since the roots have absolute value 1, inversion agrees with complex conjugation on them. Thus the coefficients of $P_v(T)$ are totally real. Hence $E$, and therefore $M$, is totally real.
		
		By \autoref{Boron}, there exist a finite CM extension $C_0/E$ and a Dirichlet density 1 set $\mathcal{L}$ of rational primes, such that for all conjugation-invariant primes $\mathfrak{p}$ of $C_0$ dividing a prime in $\mathcal{L}$, the representation $\rho_{\lambda}$ is irreducible, where $\lambda=\mathfrak{p}|_E$. Let $C$ be the Galois closure of $C_0/\mathbb{Q}$, which is also CM and contains $M$, and let $c$ be its complex conjugation.
		
		By Chebotarev, there exists a positive Dirichlet density set $\mathcal{L}'$ such that $\mathrm{Frob}_p=\{c\}$ for each $p\in\mathcal{L}'$. Since $c|_M=\mathrm{id}$, each such $p$ splits completely in $M$. Moreover, since $c$ is central in $\mathrm{Gal}(C/\mathbb{Q})$, every prime of $C$ above every $p\in\mathcal{L}'$ has the same Frobenius element $c$. Thus every prime of $C_0$ above such a $p$ is conjugation-invariant in the sense of \autoref{Boron}. Thus by this theorem, for every rational prime $p$ in $\mathcal{L}\cap\mathcal{L}'$ and every prime $\lambda$ of $E$ above $p$, the representation $\rho_{\lambda}$ is irreducible. By \autoref{Nitrogen} each such $\rho_{\lambda}$ is of type $G_2$. Since $\mathcal{L}\cap\mathcal{L}'$ has positive Dirichlet density, we have infinitely many odd rational primes $p$ satisfying (i).
	\end{proof}
	
	\begin{proof}[Proof of \autoref{Beryllium}(ii)]
		Fix $\lambda|p$. By (i), $\mathbf{G}_{\lambda}$ is of type $G_2$. Since the Galois group is compact, after a change of coordinates, $\rho_{\lambda}$ factors through $G_2(F)$ for some $p$-adic field $F$. Under these coordinates, the image $\Gamma=\mathrm{im}\rho_{\lambda}$ is a compact Zariski-dense subgroup of $G_2(F)$. For each $v\in\Sigma_K\backslash S$ with $v\nmid p$, write
		$$P_v(T)=T^7-a_vT^6+b_vT^5+\cdots,$$
		where $a_v=\mathrm{tr}(\rho_{\lambda}(\mathrm{Frob_v}))$ and $b_v=\mathrm{tr}(\wedge^2\rho_{\lambda}(\mathrm{Frob}_v))$. Since
		$$\wedge^2\rho_{\lambda}\simeq\rho_{\lambda}\oplus\mathrm{Ad}(\rho_{\lambda}),$$
		we have $\mathrm{tr}(\mathrm{Ad}\rho_{\lambda}(\mathrm{Frob}_v))=b_v-a_v$. Since $p$ splits completely in $E$, the coefficients $a_v,b_v\in\mathbb{Q}_p$ via $E\hookrightarrow E_{\lambda}\simeq\mathbb{Q}_p$. Thus $\mathrm{tr}(\mathrm{Ad}\rho_{\lambda}(\mathrm{Frob}_v))\in\mathbb{Q}_p$ for every $v\in\Sigma_K\backslash S$. By Chebotarev and continuity, $\mathrm{tr}(\mathrm{Ad}\rho_{\lambda}(\gamma))\in\mathbb{Q}_p$ for every $\gamma\in\Gamma$. Applying \autoref{Oxygen} to $\Gamma\subseteq G_2(F)$, we conclude that up to conjugacy, each $\rho_{\lambda}$ factors through $G_2(\mathbb{Q}_p)$ and has open image in it.
		
		For each $\lambda|p$, denote by $\mathfrak{g}_{\lambda}:=\mathrm{Lie}(\mathbf{G}_{\lambda})$ viewed as a $\mathrm{Gal}_K$-representation via $\mathrm{Ad}(\rho_{\lambda})$. This is irreducible since the adjoint representation of $G_2$ is irreducible and $\rho_{\lambda}(\mathrm{Gal}_K)$ is Zariski dense in $\mathbf{G}_{\lambda}$. We claim that $\mathfrak{g}_{\lambda_1}$ and $\mathfrak{g}_{\lambda_2}$ are nonisomorphic as $\mathrm{Gal}_K$-representations whenever $\lambda_1\not=\lambda_2$. Since the image of $\mathrm{Ad}(\rho_{\lambda_i})$ is Zariski dense in $\mathrm{Ad}(\mathbf{G}_{\lambda_i})(\mathbb{Q}_p)$, the algebraic monodromy group of $\mathrm{Ad}(\rho_{\lambda_i})$ is $\mathrm{Ad}(\mathbf{G}_{\lambda_i})$. Suppose for a contradiction that $\mathfrak{g}_{\lambda_1}\simeq\mathfrak{g}_{\lambda_2}$ as $\mathrm{Gal}_K$-representations. We write this isomorphism as an invertible matrix $A$. Then $\mathrm{Ad}(\mathbf{G}_{\lambda_1})\simeq\mathrm{Ad}(\mathbf{G}_{\lambda_2})$ by the conjugation by $A$. Since $G_2$ has trivial center, its adjoint map is an isomorphism onto its image. Hence the homomorphism
		$$\varphi:=\left(\xymatrix{\mathbf{G}_{\lambda_1}\ar[r]^{\mathrm{Ad}}&\mathrm{Ad}(\mathbf{G}_{\lambda_1})\ar[rr]^{x\mapsto AxA^{-1}}&&\mathrm{Ad}(\mathbf{G}_{\lambda_2})\ar[r]^{\mathrm{Ad}^{-1}}&\mathbf{G}_{\lambda_2}}\right)$$
		is an isomorphism that sends $\rho_{\lambda_1}(\gamma)$ to $\rho_{\lambda_2}(\gamma)$ for every $\gamma\in\mathrm{Gal}_K$. In particular, $\varphi\circ\rho_{\lambda_1}=\rho_{\lambda_2}$. Let $j_i:\mathbf{G}_{\lambda_i}\hookrightarrow\mathrm{GL}_7$ be the tautological representations. By uniqueness of the seven-dimensional irreducible representation of $G_2$, we obtain $j_1\simeq j_2\circ\varphi$. Hence
		$$j_1\circ\rho_{\lambda_1}\simeq j_2\circ\varphi\circ\rho_{\lambda_1}=j_2\circ\rho_{\lambda_2}.$$
		Thus $\rho_{\lambda_1}\simeq\rho_{\lambda_2}$ as $\mathrm{Gal}_K$-representations with coefficients in $\mathbb{Q}_p$. For $i=1,2$, let $\tau_i:E\hookrightarrow E_{\lambda_i}\simeq\mathbb{Q}_p$ be the embedding corresponding to $\lambda_i$. Since $\rho_{\lambda_1}\simeq\rho_{\lambda_2}$, for each $v\in\Sigma_K\backslash S$ with $v\nmid p$, we have $\tau_1(P_v(T))=\tau_2(P_v(T))$. Let $E'\subseteq E$ be the field generated by the coefficients of $P_v(T)$ for $v\in\Sigma_K\backslash S$ and $v\nmid p$. Recall that our $E$ is the fixed field of
		$$\{\sigma\in\mathrm{Aut}(\mathbb{C}):^{\sigma}\pi^{\infty}\simeq\pi^{\infty}\}.$$
		Take $\sigma\in\mathrm{Aut}(\mathbb{C}/E')$. By \cite[Th\'eor\`eme 3.13]{Cl90}, there exists a cuspidal automorphic representation $^{\sigma}\pi$ whose finite part is the scalar conjugate $^{\sigma}\pi^{\infty}$. For every $v\in\Sigma_K\backslash S$ and $v\nmid p$, $\sigma$ fixes $P_v(T)$. Hence by local-global compatibility, $^{\sigma}\pi_v\simeq\pi_v$. Since these places $v$ are cofinite, strong multiplicity one implies that $^{\sigma}\pi^{\infty}\simeq\pi^{\infty}$. Hence $\sigma$ fixes $E$. Thus $E\subseteq \mathbb{C}^{\mathrm{Aut}(\mathbb{C}/E')}=E'$. So $E=E'$. Therefore $\tau_1=\tau_2$ and hence $\lambda_1=\lambda_2$. This proves the claim.
		
		Denote by $\mathfrak{g}_p:=\mathrm{Lie}\rho_p$. Then $\mathfrak{g}_p\subseteq\oplus_{\lambda|p}\mathfrak{g}_{\lambda}$ is a $\mathrm{Gal}_K$-stable subspace under the adjoint action. The direct summands are pairwise nonisomorphic by the above claim. Since each $\rho_{\lambda}$ has open image in $G_2(\mathbb{Q}_p)$, the projection of $\mathfrak{g}_p$ to each direct summand is surjective. Hence $\mathfrak{g}_p=\oplus_{\lambda|p}\mathfrak{g}_{\lambda}$. By the $p$-adic inverse function theorem, $\rho_p$ has open image. This proves (ii).
	\end{proof}
	
	\begin{proof}[Proof of \autoref{Beryllium}(iii)]
		We construct those $v_i$ by induction. Suppose that we have already constructed $v_0,\cdots,v_{j-1}$. Fix a rational prime $p$ as in (i) with the following condition:
		$$v_{\mathfrak{p}}\left(\sigma(\Delta_{v_i})\right)=0,\ \forall\mathfrak{p}|p\ \text{in}\ M,\ \sigma\in X,\ 0\le i\le j-1.$$
		Such a prime exists since this only rules out finitely many primes. For integers $m\ge1$ and $\delta\in\{0,1\}$, define
		$$h_{m,\delta}=\left(\begin{array}{cc}1&p^m\\p^{m+\delta}&1+p^{2m+\delta}\end{array}\right)\in\mathrm{SL}_2(\mathbb{Q}_p).$$
		Let $\iota:\mathrm{SL}_2\to G_2$ be the short-root homomorphism. Write $z$ and $z^{-1}$ as the eigenvalues of $h_{m,\delta}$. Then the eigenvalues of $\iota(h_{m,\delta})$ are
		$$z,z^{-1},z,z^{-1},z^2,z^{-2},1.$$
		Since $z+z^{-1}=\mathrm{tr}(h_{m,\delta})=2+p^{2m+\delta}$, we have
		$$\Delta(\iota(h_{m,\delta}))=(4+p^{2m+\delta})p^{3(2m+\delta)}.$$
		As $p$ is odd, $4+p^{2m+\delta}$ is a unit in $\mathbb{Z}_p$. Hence
		$$v_p\left(\Delta(\iota(h_{m,\delta}))\right)=3(2m+\delta).$$
		Fix a prime $\lambda_0|p$ of $E$. Take $a_m=(a_{\lambda})_{\lambda|p}\in\prod_{\lambda|p}G_2(\mathbb{Q}_p)$ such that $a_{\lambda_0}=\iota(h_{m,1})$ and $a_{\lambda}=\iota(h_{m,0})$ for each $\lambda\not=\lambda_0$. This tuple $a_m$ tends to the identity as $m$ tends to infinity. By part (ii), we take $m$ sufficiently large so that $a_m\in\mathrm{im}\rho_p$. Continuity of $\Delta$ and the valuation topology of $\mathbb{Q}_p$ imply that
		$$U=\left\{(g_{\lambda})_{\lambda|p}:v_p(\Delta(g_{\lambda}))=v_p(\Delta(a_{\lambda}))\ \text{for every}\ \lambda|p\right\}$$
		is nonempty and open in $\prod_{\lambda|p}G_2(\mathbb{Q}_p)$. Since taking conjugation does not change characteristic polynomials, this $U$ is conjugation-invariant. Thus $\rho_p^{-1}(U)$ is a nonempty open conjugation-invariant subset of $\mathrm{Gal}_K$. By Chebotarev, there exists $v_j\in\Sigma_K\backslash S$ with $v_j\nmid p$ such that $\rho_p(\mathrm{Frob}_{v_j})\in U$. Consequently,
		$$v_{\lambda}\left(\Delta_{v_j}\right)=\begin{cases}6m+3&\lambda=\lambda_0\\6m&\lambda\not=\lambda_0\end{cases}$$
		We now check that $[\sigma(\Delta_{v_i})]$ for $0\le i\le j$ and $\sigma\in X$ are linearly independent over $\mathbb{F}_2$ in $M^{\times}/(M^{\times})^2$. Suppose that we have integers $a_{i,\sigma}\in\{0,1\}$ for each $0\le i\le j$ and $\sigma\in X$, such that
		\begin{equation}\label{2}
		\prod_{0\le i\le j}\prod_{\sigma\in X}\sigma(\Delta_{v_i})^{a_{i,\sigma}}\in (M^{\times})^2.
		\end{equation}
		Fix a $\sigma'\in X$ and take a prime $\mathfrak{p}$ of $M$ such that $(\sigma')^{-1}(\mathfrak{p})=\lambda_0$. For this fixed $\mathfrak{p}$, the map $\sigma\mapsto\sigma^{-1}(\mathfrak{p})$ gives a one-to-one correspondence from $X$ to the primes $\lambda$ of $E$ above $p$, since $p$ splits completely in $M$. Moreover, $v_{\mathfrak{p}}(\sigma(x))=v_{\sigma^{-1}(\mathfrak{p})}(x)$ for $x\in E^{\times}$. Hence
		$$v_{\mathfrak{p}}\left(\sigma(\Delta_{v_i})\right)=\begin{cases}0&0\le i\le j-1\\6m+3&i=j,\sigma=\sigma'\\6m&i=j,\sigma\not=\sigma'\end{cases}$$
		Hence applying $v_{\mathfrak{p}}(-)$ to (\ref{2}) yields
		$$(6m+3)a_{j,\sigma'}+6m\sum_{\sigma\not=\sigma'}a_{j,\sigma}\in2\mathbb{Z}.$$
		Hence $a_{j,\sigma'}=0$. Since we can pick any $\sigma'$ in $X$, the exponent $a_{j,\sigma}=0$ for all $\sigma\in X$. Hence from (\ref{2}) we have
		$$\prod_{0\le i\le j-1}\prod_{\sigma\in X}\sigma(\Delta_{v_i})^{a_{i,\sigma}}\in (M^{\times})^2.$$
		Then as $[\sigma(\Delta_{v_i})]$ are linearly independent for $0\le i\le j-1$ and $\sigma\in X$, we have $a_{i,\sigma}=0$ for them. This completes the proof.
	\end{proof}
	
	Let $\mathcal{L}$ be the set of rational primes $\ell$ for which $\rho_{\lambda}$ admits a decomposition of dimensional type $3+3+1$ for at least one prime $\lambda$ of $E$ above $\ell$. The goal is to show $\mathcal{L}$ has Dirichlet density zero. Fix a positive integer $N$, take $(N+1)$ primes $v_0,\cdots,v_N$ of $K$ and put $r_i=\Delta_{v_i}/\Delta_{v_0}$ as in \autoref{Beryllium}(iii). Define the number field
	$$L=M\left(\left\{\sqrt{\sigma(r_i)}:\sigma\in X,1\le i\le N\right\}\right).$$
	Since $M/\mathbb{Q}$ is Galois and the square classes defining $L/M$ are stable under $\mathrm{Gal}(M/\mathbb{Q})$, Kummer theory implies that $L/\mathbb{Q}$ is Galois. For each $g\in\mathrm{Gal}(L/\mathbb{Q})$, the restriction $g|_M$ acts on the set $X=\mathrm{Hom}_{\mathbb{Q}}(E,M)$. For each $\sigma\in X$, denote by $f_{g,\sigma}$ the size of the $g|_M$-orbit of $X$ that contains $\sigma$. Define a subset of $\mathrm{Gal}(L/\mathbb{Q})$:
	$$C_N=\left\{g\in\mathrm{Gal}(L/\mathbb{Q}):\exists\sigma\in X\ \text{such that}\ g^{f_{g,\sigma}}(\sqrt{\sigma(r_i)})=\sqrt{\sigma(r_i)}\ \text{for all}\ 1\le i\le N\right\}.$$
	\begin{lem}\label{Fluorine}
		Remove the rational primes ramified in $L$ and the residue characteristics of $v_0,\cdots,v_N$ from $\mathcal{L}$. Then for every $N\ge1$, $C_N$ is conjugation-invariant and contains the conjugacy class $\mathrm{Frob}_{\ell}$ for each $\ell\in\mathcal{L}$.
	\end{lem}
	\begin{proof}
		Fix $N$ and write $C=C_N$. Take $\ell\in\mathcal{L}$ as stated, and pick a prime $\lambda|\ell$ of $E$ such that $\rho_{\lambda}$ decomposes as $3+3+1$. We first show that each $r_i$ is a square in $E_{\lambda}$ for $1\le i\le N$. Take $\gamma\in\mathrm{Gal}(\overline{E}_{\lambda}/E_{\lambda})$. For each $v\in\Sigma_K\backslash S$ with $v\nmid\ell$, the characteristic polynomials of $\rho_{\lambda}(\mathrm{Frob}_v)$ and $^{\gamma}\rho_{\lambda}(\mathrm{Frob}_v)$ are equal since $P_v(T)$ has coefficients in $E\subseteq E_{\lambda}$. Hence Chebotarev and Brauer-Nesbitt imply $\rho_{\lambda}\simeq^{\gamma}\rho_{\lambda}$. Thus $\gamma$ either swaps $U_{\lambda}$ and $U_{\lambda}^{\vee}$ or fixes them up to isomorphism. Consequently, $\gamma$ either swaps the traces $a_{\lambda,v_i}$ and $b_{\lambda,v_i}$ for all $0\le i\le N$, or fixes them for all $i$. This implies that $\gamma$ fixes the ratios
		$$\frac{a_{\lambda,v_i}-b_{\lambda,v_i}}{a_{\lambda,v_0}-b_{\lambda,v_0}}$$
		for each $1\le i\le N$. We conclude that these ratios lie in $\overline{E}_{\lambda}^{\mathrm{Gal}(\overline{E}_{\lambda}/E_{\lambda})}=E_{\lambda}$. Hence each $r_i$ is a square in $E_{\lambda}$ for $1\le i\le N$.
		
		Fix a prime $\mathfrak{p}$ of $M$ above our fixed $\ell$. Take $\sigma\in X$ such that $\sigma^{-1}(\mathfrak{p})=\lambda$. Put $h=\mathrm{Frob}_{\mathfrak{p}}\in\mathrm{Gal}(M/\mathbb{Q})$. Then $h$ acts on $X$ via composition. The $h$-orbit containing $\sigma$ is the cycle
		$$\{\sigma,h\circ\sigma\cdots,h^{f-1}\circ\sigma\},$$
		where $f$ is the residue degree of $\lambda$ over $\ell$. Choose a prime $\mathfrak{P}$ of $L$ above $\mathfrak{p}$, and put $g=\mathrm{Frob}_{\mathfrak{P}}\in\mathrm{Gal}(L/\mathbb{Q})$. Since $h^f\circ\sigma=\sigma$, the restriction $h^f|_{\sigma(E)}=\mathrm{id}$. Hence $g^f|_{\sigma(E)}=\mathrm{id}$. Thus $g^f\in\mathrm{Gal}(L/\sigma(E))$. The embedding $\sigma:E\hookrightarrow M$ extends to the completions $E_{\lambda}\hookrightarrow M_{\mathfrak{p}}$. $g^f$ is the Frobenius element at $\mathfrak{P}$ for $L/\sigma(E)$. It therefore acts on $L_{\mathfrak{P}}$ over $\sigma(E_{\lambda})$. Hence $g^f|_{\sigma(E_{\lambda})}=\mathrm{id}$. As $r_i$ is a square in $E_{\lambda}$, both square roots of $\sigma(r_i)$ lie in $\sigma(E_{\lambda})$. Hence $g^f$ fixes each $\sqrt{\sigma(r_i)},1\le i\le N$ for our chosen $\sigma$. Thus, for the stated $\ell$, we have a Frobenius element $g$ in the conjugacy class $\mathrm{Frob}_{\ell}$ such that $g\in C$.
		
		We show that $C$ is conjugation-invariant. Take $g\in C$. Then there is an embedding $\sigma\in X$ such that $g^f$ fixes all $\sqrt{\sigma(r_i)}$ for $1\le i\le N$, where $f$ is the length of the $g|_M$-orbit of $X$ containing $\sigma$. For any $t\in\mathrm{Gal}(L/\mathbb{Q})$, take the embedding $\sigma':=t|_M\circ\sigma$. Then the $(tgt^{-1})|_M$-orbit of $X$ containing $\sigma'$ again has length $f$, since composition by $t|_M$ gives a bijection from the $g|_M$-orbit of $\sigma$ to the $(tgt^{-1})|_M$-orbit of $\sigma'$. Since $g^f$ fixes all $\sqrt{\sigma(r_i)}$, the conjugate $(tgt^{-1})^f$ fixes each $\sqrt{\sigma'(r_i)}$. Hence $tgt^{-1}\in C$. It follows that the whole conjugacy class $\mathrm{Frob}_{\ell}$ is contained in $C$. This proves the lemma.
	\end{proof}
	
	\begin{lem}\label{Neon}
		For $h\in\mathrm{Gal}(M/\mathbb{Q})$, denote by $c(h)$ the number of orbits of $h$ acting on $X$. Then
		$$\frac{|C_N|}{|\mathrm{Gal}(L/\mathbb{Q})|}=\frac{1}{[M:\mathbb{Q}]}\sum_{h\in\mathrm{Gal}(M/\mathbb{Q})}\left(1-(1-2^{-N})^{c(h)}\right)\le\frac{[E:\mathbb{Q}]}{2^N}$$
	\end{lem}
	\begin{proof}
		Write $C=C_N$. For each $h\in\mathrm{Gal}(M/\mathbb{Q})$, let $C_h$ consist of $g\in C$ such that $g|_M$ is $h$. We now calculate the size $|C_h|$. Fix a lift $g_0$ of $h$ to $\mathrm{Gal}(L/\mathbb{Q})$. Then each lift $g$ of $h$ has a unique form $g=ug_0$ with $u\in\mathrm{Gal}(L/M)$. For each $\sigma\in X$ and each $1\le i\le N$, choose a square root $x_{\sigma,i}=\sqrt{\sigma(r_i)}$ and define $\varepsilon_{\sigma,i}(u)\in\{\pm1\}$ by $u(x_{\sigma,i})=\varepsilon_{\sigma,i}(u)x_{\sigma,i}$. By \autoref{Beryllium}(iii), the classes of the $\sigma(r_i)$ are linearly independent in $M^{\times}/(M^{\times})^2$. Thus, by Kummer theory, as $u$ varies, every assignment of these signs occurs exactly once. 
		
		For each $h$-orbit $\mathcal{O}$ in $X$, fix an element $\tau\in\mathcal{O}$, put $f=|\mathcal{O}|$, and write $g_0^f(x_{\tau,i})=\delta_{\mathcal{O},i}x_{\tau,i}$ where $\delta_{\mathcal{O},i}\in\{\pm1\}$. Since $g_0$ sends $x_{\sigma,i}$ to $\pm x_{h\circ\sigma,i}$, iterating $g=ug_0$ gives
		$$g^f(x_{\tau,i})=\delta_{\mathcal{O},i}\left(\prod_{\sigma\in\mathcal{O}}\varepsilon_{\sigma,i}(u)\right)x_{\tau,i}$$
		Since $g_0^f$ commutes with $g_0$, the product and $\delta_{\mathcal{O},i}$ are independent of the representative $\tau$. Thus the lifts $g$ to $h$ for which $g^f$ fixes each $\sqrt{\tau(r_i)},1\le i\le N$ correspond to sign assignments of $\varepsilon_{\sigma,i}(u)\in\{\pm1\}$ for $\sigma\in\mathcal{O},1\le i\le N$ satisfying
		$$\prod_{\sigma\in\mathcal{O}}\varepsilon_{\sigma,i}(u)=\delta_{\mathcal{O},i},\quad\left(1\le i\le N\right).$$
		For each $i$, all but one of the signs $\varepsilon_{\sigma,i}(u)$ with $\sigma\in\mathcal{O}$ can be chosen freely, and the product equation determines the remaining sign. Thus each equation has $2^{|\mathcal{O}|-1}$ solutions. Since the $N$ equations involve disjoint sets of variables, there are $2^{N(|\mathcal{O}|-1)}$ assignments on $\mathcal{O}$ satisfying all of them. The signs outside $\mathcal{O}$ are unrestricted. Thus the number $d_{\mathcal{O}}$ of such assignments is
		$$d_{\mathcal{O}}=2^{N(|\mathcal{O}|-1)}\cdot2^{N([E:\mathbb{Q}]-|\mathcal{O}|)}=2^{N([E:\mathbb{Q}]-1)}.$$
		The total number of assignments of signs is $[L:M]=2^{N[E:\mathbb{Q}]}$. Hence $1-d_{\mathcal{O}}/[L:M]$ is the fraction of all lifts $g$ of $h$ for which $g^{|\mathcal{O}|}$ fails to fix at least one of these square roots. Since the systems of equations for different orbits have disjoint variables, the fraction of lifts failing to fix some square root on every orbit is $\prod_{\mathcal{O}}(1-d_{\mathcal{O}}/[L:M])$. We conclude that
		$$|C_h|=[L:M]\left(1-\prod_{\mathcal{O}}\left(1-\frac{d_{\mathcal{O}}}{[L:M]}\right)\right)=[L:M]\left(1-(1-2^{-N})^{c(h)}\right).$$
		The number of $h$-orbits $c(h)$ can be at most $|X|=[E:\mathbb{Q}]$. Applying Bernoulli's inequality $(1+t)^a\ge1+at$ for $t\ge-1$ and $a\ge1$ to $t=-2^{-N}$ and $a=c(h)$, we obtain
		$$\left(1-(1-2^{-N})^{c(h)}\right)\le\left(1-(1-c(h)2^{-N})\right)=c(h)2^{-N}\le[E:\mathbb{Q}]2^{-N}.$$
		Hence
		\begin{align*}\frac{|C|}{[L:\mathbb{Q}]}&=\sum_{h\in\mathrm{Gal}(M/\mathbb{Q})}\frac{|C_h|}{[L:\mathbb{Q}]}\\&=\frac{1}{[M:\mathbb{Q}]}\sum_{h\in\mathrm{Gal}(M/\mathbb{Q})}\left(1-(1-2^{-N})^{c(h)}\right)\\&\le\frac{1}{[M:\mathbb{Q}]}\sum_{h\in\mathrm{Gal}(M/\mathbb{Q})}\left([E:\mathbb{Q}]2^{-N}\right)\\&=\frac{[E:\mathbb{Q}]}{2^N}\end{align*}
	\end{proof}
	
	By \autoref{Fluorine}, for each $N\ge1$, the Frobenius conjugacy class of every $\ell\in\mathcal{L}$ outside a finite set is contained in the conjugation-invariant set $C_N$.  Hence Chebotarev and \autoref{Neon} yield
	$$\overline{\delta}(\mathcal{L})\le\frac{|C_N|}{[L:\mathbb{Q}]}\le\frac{[E:\mathbb{Q}]}{2^N},$$
	where $\overline{\delta}$ denotes the upper Dirichlet density. As $N$ can be any positive integer, $\mathcal{L}$ has Dirichlet density zero. This finishes the proof of \autoref{main}.
	
	\section{Extension of the Feng-Whitmore Theorem}
	
	In this section, we prove \autoref{main2} by following \cite{FW25} and treating the additional cases that arise when the assumption $7\nmid n$ is removed. Let $(\pi,\chi)$ and $\{\rho_{\lambda}\}_{\lambda}$ as in \autoref{main2}. Throughout this section, we retain the assumption
	\begin{equation}\label{6}
	4|n\quad\Rightarrow\quad n=4\text{ or }n=4p\text{ for some prime }p.
	\tag{N}
	\end{equation}
	
	\subsection{Reduction to simple derived monodromy}
	
	Following the reduction step of \cite{Xi19}, as reformulated in \cite[Propositions 4.13 and 4.14]{Hu23b}, we may reduce the proof of \autoref{main2} to the case where $K$ is imaginary CM and is the maximal CM subfield of the connectedness field given by \autoref{Hydrogen}(i) for the auxiliary compatible system
	$$\left\{\mathrm{Ind}_K^{K^+}\rho_{\lambda}\oplus\mu_{\lambda}\right\}_{\lambda}.$$
	The reduced system has dimension dividing $n$, so the numerical condition (\ref{6}) is preserved under this reduction. By \cite[Corollary 1]{Xi19} (see also \cite[Theorem 3.4]{FW25}), there exists a set $\mathcal{L}'$ of rational primes of positive Dirichlet density such that $\rho_{\lambda}$ is Lie-irreducible for every $\ell\in\mathcal{L}'$ and every $\lambda|\ell$.
	
	By \cite[Proposition 4.5]{FW25}, for sufficiently large $\ell\in\mathcal{L}'$ and every $\lambda|\ell$, we have a tensor decomposition
	$$\rho_{\lambda}\simeq\bigotimes_{1\le i\le s}\rho_{i,\lambda},$$
	where each $\rho_{i,\lambda}$ is Lie-irreducible and polarized, with simple derived monodromy Lie algebra $\mathfrak{g}_{i,\lambda}$, such that
	$$\mathrm{Lie}(\mathbf{G}_{\lambda}^{\mathrm{der}})\simeq\prod_{1\le i\le s}\mathfrak{g}_{i,\lambda}.$$
	
	After replacing $\mathcal{L}'$ by a positive Dirichlet density subset, \cite[Corollary 4.6]{FW25} implies that, for every $\ell\in\mathcal{L}'$ and every $\lambda|\ell$, all the tensor factors of $\rho_{\lambda}$ become automorphic after restriction to a common finite imaginary CM extension $K'/K$, which may depend on $\lambda$. Its proof applies without the assumption that $7\nmid n$. Indeed, after shrinking $\mathcal{L}'$ as in their proof, every tensor factor of odd dimension or dimension two is totally odd by \cite[Lemma 2.9(1) and (2)]{FW25}. Since $n=\prod_{1\le i\le s}\dim\rho_{i,\lambda}$, condition (\ref{6}) implies that at most one tensor factor has even dimension greater than two. As $\rho_{\lambda}$ is totally odd, \cite[Lemma 2.9(3)]{FW25} implies that all the tensor factors are totally odd. Hence the same potential automorphy argument applies.
	
	Since the rank of the derived monodromy Lie algebra of $\rho_{\lambda}$ is independent of $\lambda$, we may choose an infinite set $\Lambda$ of coefficient primes lying above primes in $\mathcal{L}'$, whose underlying rational primes have positive upper Dirichlet density, such that the number $s$ of simple factors of the derived monodromy Lie algebra $\mathfrak{g}_{\lambda}$ of $\rho_{\lambda}$ is independent of $\lambda\in\Lambda$.
	
	Fix a coefficient prime $\lambda'\in\Lambda$. After applying potential automorphy to the tensor factors of $\rho_{\lambda'}$ and enlarging the coefficient field to a common imaginary CM field $E$, we thus obtain compatible systems $\{r_{i,\lambda}\}_{\lambda}$ attached to regular algebraic polarized cuspidal automorphic representations over $K'$, such that
	$$\rho_{\lambda}|_{\mathrm{Gal}_{K'}}\simeq\bigotimes_{1\le i\le s}r_{i,\lambda}$$
	for every $\lambda$. We replace each place in $\Lambda$ by a place of the enlarged coefficient field above it. Then denote by $\mathfrak{g}_{i,\lambda}$ the derived monodromy Lie algebra of $r_{i,\lambda}$. The tensor decomposition induces an embedding
	$$\mathfrak{g}_{\lambda}\hookrightarrow\prod_{1\le i\le s}\mathfrak{g}_{i,\lambda}$$
	whose projections to each factor are surjective. By construction, this embedding is an isomorphism at the coefficient place $\lambda'$ used to obtain the compatible systems $\{r_{i,\lambda}\}_{\lambda}$. Hence \autoref{Hydrogen}(ii) implies that
	$$\rk\mathfrak{g}_{\lambda}=\sum_{1\le i\le s}\rk\mathfrak{g}_{i,\lambda}$$
	for every $\lambda$. By \cite[Lemma 4.9]{FW25}, we obtain
	$$\mathfrak{g}_{\lambda}\simeq\prod_{1\le i\le s}\mathfrak{g}_{i,\lambda}.$$
	By $\lambda$-independence of rank, each $\mathfrak{g}_{i,\lambda}$ has positive rank, since $\mathfrak{g}_{i,\lambda'}$ is simple. Since $\mathfrak{g}_{\lambda}$ has exactly $s$ simple factors for every $\lambda\in\Lambda$, each $\mathfrak{g}_{i,\lambda}$ is simple for these $\lambda$. Moreover, each $r_{i,\lambda}$ is Lie-irreducible for $\lambda\in\Lambda$, since $\rho_{\lambda}$ is Lie-irreducible. By \cite[Proposition 5.3.2]{BLGGT14} and the positive upper Dirichlet density of the rational primes underlying $\Lambda$, we may therefore choose a sufficiently large $\lambda_0\in\Lambda$ at which all $r_{i,\lambda_0}$ satisfy the hypotheses of \cite[Theorem 4.8]{FW25}, without the condition $7\nmid\dim r_{i,\lambda_0}$. The intersection of the finitely many Dirichlet density-one sets on which the systems $\{r_{i,\lambda}\}_{\lambda}$ are Lie-irreducible is a Dirichlet density-one set on which $\rho_{\lambda}$ is Lie-irreducible. This reduces the proof of \autoref{main2} to the following \cite[Theorem 4.8]{FW25} without the condition $7\nmid n$:
	\begin{thm}[\cite{FW25}, Theorem 4.8, with the condition $7\nmid n$ removed]\label{Phosphorus}
		Let $K$ be an imaginary CM field and let $(\pi,\chi)$ be a regular algebraic polarized cuspidal automorphic representation of $\mathrm{GL}_n(\mathbb{A}_K)$ with associated $E$-rational compatible system $\{\rho_\lambda\}_\lambda$, where $E$ is an imaginary CM field. Assume that if $4\mid n$, then $n=4p$ for some prime $p$. Suppose that there exists a coefficient prime $\lambda_0|\ell_0$ satisfying the following conditions:
		\begin{enumerate}
			\item[(i)] $\rho_{\lambda_0}$ is Lie-irreducible.
			\item[(ii)] The derived monodromy Lie algebra $\mathfrak{g}_{\lambda_0}$ is nonzero and simple.
			\item[(iii)] $\ell_0\geq 2n+2$ and $\ell_0$ is unramified in $K$.
			\item[(iv)] $\rho_{\lambda_0}$ is crystalline at all places above $\ell_0$, and for each embedding $\tau:K\hookrightarrow\overline{\mathbb{Q}}_{\ell_0}$, there exists an integer $a_\tau$ such that
			$$\mathrm{HT}_{\tau}(\rho_{\lambda_0})\subseteq [a_\tau,a_\tau+2\ell_0-4].$$
			\item[(v)] The restriction $\overline{\rho}_{\lambda_0}|_{\mathrm{Gal}_{K(\zeta_{\ell_0})}}$ is irreducible.
		\end{enumerate}
		Then there exists a set $\mathcal{L}$ of rational primes of Dirichlet density one such that $\rho_\lambda$ is Lie-irreducible for every $\ell\in\mathcal{L}$ and every $\lambda|\ell$.
	\end{thm}
	
	\subsection{The simple Lie algebra cases}
	
	By the preceding reduction, it suffices to prove \autoref{Phosphorus}. Let $\lambda_0$ be as in that theorem, and denote by $\mathfrak{g}=\mathrm{Lie}(\mathbf{G}_{\lambda_0}^{\mathrm{der}})$ the simple derived monodromy Lie algebra of $\rho_{\lambda_0}$, with its tautological representation $V$. 
	
	Since $\rho_{\lambda_0}$ is Hodge-Tate regular and Lie-irreducible, \cite[Lemma 2.20]{FW25} implies that the formal character of $V$ is multiplicity-free. By the classification of irreducible multiplicity-free representations of simple Lie algebras \cite[Theorem 4.6.3]{Ho95}, the pair $(\mathfrak{g},V)$ belongs to the following list, up to isomorphism and duality:
	\begin{center}
		\begin{tabular}{c|l}
			$\mathfrak{g}$ & $V$ (up to duality) \\
			\hline
			$A_m$ & $\operatorname{Sym}^a\mathrm{Std},\ \bigwedge^b\mathrm{Std}$ \\
			$B_m$ ($m\geq 2$) & Standard or spin representation \\
			$C_m$ ($m\geq 3$) & Standard representation; $V_{\omega_3}$ if $m=3$ \\
			$D_m$ ($m\geq 4$) & Standard or half-spin representations \\
			$E_6$ & $V_{\omega_1}, V_{\omega_6}$ (dimension $27$) \\
			$E_7$ & $V_{\omega_7}$ (dimension $56$) \\
			$G_2$ & $V_{\omega_1}$ (dimension $7$) \\
		\end{tabular}
	\end{center}
	
	We now summarize the proofs of the corresponding cases in \cite{FW25}.
	\begin{itemize}
		\item \textbf{Type $A$}. Propositions 5.14-5.17 treat these cases. The case of $\mathfrak{sl}_4$ acting on $\bigwedge^2\mathrm{Std}$ is covered by \cite[Proposition 5.20]{FW25}, since $\mathfrak{sl}_4\simeq\mathrm{so}_6$ and the exterior square corresponds the standard representation of $\mathfrak{so}_6$. The proofs of Propositions 5.16 and 5.17 use the condition $7\nmid n$ to exclude the proper exterior powers of the standard representation of $\mathfrak{sl}_8$. Proposition 5.16 also uses this condition to exclude a possible derived monodromy Lie algebra of type $\mathfrak{e}_7$ at another coefficient place. This possibility is already excluded by (\ref{6}), since $n=56$. The remaining exterior powers have dimensions 28 and 70.
		\item \textbf{Types} $B,C,D$. The standard representations are treated in Proposition 5.20. Condition (\ref{6}) excludes spin and half-spin representations of dimension at least 16. The remaining spin case of type $B_3$ is treated in Proposition 5.22, and the eight-dimensional half-spin representations of $D_4$ are treated by Proposition 5.20 using triality \cite[\S 20.3]{FH13}. So only the fourteen-dimensional representation of $C_3$ remains.
		\item \textbf{Exceptional types}. The case of $\mathfrak{e}_6$ is treated in Proposition 5.23. The case of $\mathfrak{e}_7$ is excluded by (\ref{6}). So the remaining case is the seven-dimensional standard representation of $\mathfrak{g}_2$.
	\end{itemize}
	
	Consequently, to complete the proof of \autoref{main2}, it remains to prove the following proposition. These cases are treated without the additional local and residual irreducibility conditions imposed in \cite[Theorem 4.8]{FW25}.
	\begin{prop}\label{Silicon}
		Let $K$ be imaginary CM and let $(\pi,\chi)$ be a regular algebraic polarized cuspidal automorphic representation of $\mathrm{GL}_n(\mathbb{A}_K)$ with associated $E$-rational compatible system $\{\rho_{\lambda}\}_{\lambda}$. Suppose that there exists a coefficient prime $\lambda_0$ such that $\rho_{\lambda_0}$ is Lie-irreducible, and the Lie algebra $\mathfrak{g}$ of $\mathbf{G}_{\lambda_0}^{\mathrm{der}}$ with its tautological representation is one of the following:
		\begin{enumerate}
			\item $\mathfrak{g}=\mathfrak{g}_2$ and $V=V_{\omega_1}$ of dimension seven.
			\item $\mathfrak{g}=\mathfrak{sp}_6$ and $V=V_{\omega_3}$ of dimension fourteen.
			\item $\mathfrak{g}=\mathfrak{sl}_8$ and $V$ is either $\bigwedge^2\mathrm{Std}$, its dual, or $\bigwedge^4\mathrm{Std}$.
		\end{enumerate}
		Then there exists a set $\mathcal{L}$ of rational primes of Dirichlet density one such that $\rho_{\lambda}$ is Lie-irreducible for every $\ell\in\mathcal{L}$ and every $\lambda|\ell$.
	\end{prop}
	
	Case (i) follows from our main theorem \autoref{main}. We treat cases (ii) and (iii) in the following two subsections.
	
	\subsection{The fourteen-dimensional representation of type $C_3$}
	
	In this subsection, we prove \autoref{Silicon}(ii). We first use the geometry of the derived formal character to classify the possible decompositions of $\rho_{\lambda}^{\mathrm{der}}$.
	\begin{lem}
		If the tautological representation of $\mathbf{G}_{\lambda_0}^{\mathrm{der}}$ is the unique multiplicity-free fourteen-dimensional irreducible representation $V_{\omega_3}$ of $\mathfrak{sp}_6$, then for each $\lambda$, the tautological representation of $\mathbf{G}_{\lambda}^{\mathrm{der}}$ is one of the following.
		\begin{enumerate}
			\item[(i)] $\left(\mathrm{Sp}_6, V_{\omega_3}\right)$.
			\item[(ii)] $\mathrm{SL}_2\times\mathrm{SL}_2\times\mathrm{SL}_2$ acting by
			$$\left(\mathrm{Std}\boxtimes\mathbf{1}\boxtimes\mathbf{1}\right)\oplus\left(\mathbf{1}\boxtimes\mathrm{Std}\boxtimes\mathbf{1}\right)\oplus\left(\mathbf{1}\boxtimes\mathbf{1}\boxtimes\mathrm{Std}\right)\oplus\left(\mathrm{Std}\boxtimes\mathrm{Std}\boxtimes\mathrm{Std}\right).$$
			\item[(iii)] $\mathrm{Sp}_4\times\mathrm{SL}_2$ acting by
			$$\left(\mathrm{Std}\boxtimes\mathbf{1}\right)\oplus\left(W\boxtimes\mathrm{Std}\right),$$ 
			where $W$ is the unique five-dimensional irreducible summand of $\bigwedge^2(\mathrm{Std})$ of $\mathrm{Sp}_4$.
		\end{enumerate}
	\end{lem}
	\begin{proof}
		The formal character of $\mathbf{G}_{\lambda_0}^{\mathrm{der}}$ is
		$$\{\pm x,\pm y,\pm z,\pm x\pm y\pm z\},$$
		where $x$, $y$ and $z$ are linearly independent characters. By \autoref{Hydrogen}, all $\mathbf{G}_{\lambda}^{\mathrm{der}}$ have the same formal character as above. With respect to the basis $x,y,z$, the weights $\pm x\pm y\pm z$ are the eight vertices of the cube $[-1,1]^3$, and the weights $\pm x, \pm y, \pm z$ are the centers of its six faces. Each Weyl reflection permutes the weights and hence preserves their convex hull $[-1,1]^3$. It therefore permutes the eight vertices and the six face centers separately. Since it permutes $\{\pm x,\pm y,\pm z\}$, it either changes the sign of one coordinate and fixes the other two coordinates, or swaps two coordinates with equal signs and fixes the third. Thus the direction of this reflection is either a coordinate axis or a diagonal line in a coordinate plane. As each root is a difference of two weights, the possible roots are:
		$$\begin{array}{c|l}
			\text{Length} & \text{Candidate roots} \\ \hline
			\sqrt{2}
			& \pm x\pm y,\ \pm y\pm z,\ \pm z\pm x \\[2pt]
			2
			& \pm 2x,\ \pm 2y,\ \pm 2z \\[2pt]
			2\sqrt{2}
			& \pm 2x\pm 2y,\ \pm 2y\pm 2z,\ \pm 2z\pm 2x
		\end{array}$$
		
		Equip the weight space with the inner product for which $x,y,z$ are perpendicular to each other and each has length one. The length-$2\sqrt{2}$ candidates cannot be roots: by symmetry of the weight set, it suffices to check $2x+2y$:
		$$\langle x,(2x+2y)^{\vee}\rangle=\frac{2(x,2x+2y)}{(2x+2y,2x+2y)}=\frac{1}{2},$$
		which is not an integer. Hence they cannot be roots.
		
		We refer to the length-$\sqrt{2}$ candidates as diagonal candidates and the length-$2$ candidates as coordinate candidates. We now classify the possible root systems of $\mathbf{G}_{\lambda}^{\mathrm{der}}$.
		
		\noindent\textbf{Claim.} If both diagonal lines in a
			coordinate plane contain roots, then both coordinate axes in that plane also contain roots. Hence in this coordinate plane the roots are the $C_2$ root system.
			
		After relabeling the coordinates, suppose $\alpha=x+y$ and $\beta=x-y$ are roots. We want to show that both $2x$ and $2y$ are roots. Suppose for a contradiction that neither $2x$ nor $2y$ is a root. For each root $\gamma$, let $H_{\gamma}$ be its coroot element and choose nonzero $E_{\gamma}\in\mathfrak{g}_{\gamma}$ and $F_{\gamma}\in\mathfrak{g}_{-\gamma}$ in the root spaces such that $[E_{\gamma},F_{\gamma}]=H_{\gamma}$. Since $\alpha-\beta=2y$ is not a root, we have $[E_{\alpha},F_{\beta}]=0$.
		
		Choose nonzero vectors $a\in V_{z+\alpha}$ and $b\in V_{z+\beta}$. The weights $z+2\alpha$ and $z+2\beta$ are absent, so $E_{\alpha}a=E_{\beta}b=0$. We have
		$$H_{\alpha}a=(z+\alpha)(H_{\alpha})a=\left(\langle z,\alpha^{\vee}\rangle+\langle\alpha,\alpha^{\vee}\rangle\right)a=(0+2)a=2a.$$
		Similarly, $H_{\beta}b=2b$. Thus
		$$E_{\alpha}(F_{\alpha}a)=2a,\quad E_{\beta}(F_{\beta}b)=2b.$$
		In particular, $F_{\alpha}a$ and $F_{\beta}b$ are two nonzero vectors in the weight space $V_z$. Since the weight $z+\alpha+\beta=z+2x$ is absent, $E_{\alpha}b=0$. Hence $E_{\alpha}(F_{\beta}b)=F_{\beta}(E_{\alpha}b)=0$ while $E_{\alpha}(F_{\alpha}a)\not=0$. So the vectors $F_{\alpha}a$ and $F_{\beta}b$ are linearly independent, contradicting that $\dim V_z=1$. Thus at least one of $2x$ and $2y$ occurs as a root. The Weyl reflection $s_{x-y}$ exchanges them, so both are roots. This proves the claim. 
		
		We now distinguish three cases.
		\begin{itemize}
			\item No diagonal candidates occur as roots. Since $\mathbf{G}_{\lambda}^{\mathrm{der}}$ has rank three, the root system consists of all six coordinate candidates. This is $A_1^3$. Choose the simple roots as $2x,2y,2z$. The highest weights in the formal character are $x,y,z,x+y+z$. The first three correspond to the standard representation of each $A_1$ factor with trivial action on the other two factors. The last highest weight corresponds to the external tensor product of the three standard representations of each $A_1$ factor. This gives case (ii).
			\item Exactly one coordinate plane contains diagonal roots. After relabeling the coordinates, we may assume this coordinate plane is the $xy$-plane. Since $\mathbf{G}_{\lambda}^{\mathrm{der}}$ has rank three, $2z$ must be a root and the roots in the $xy$-plane span that plane. In the $xy$-plane, either we have two linearly independent diagonal roots, which give the $C_2$ root system by the claim above, or we have one diagonal root and one coordinate root, whose reflections give the other coordinate axis and the other diagonal line, again giving $C_2$ root system. In either case, the root system is of type $C_2\times A_1$. Choose the simple roots as $x-y,2y,2z$. Then the highest weights in the formal character of $\mathbf{G}_{\lambda}^{\mathrm{der}}$ are $x$ and $x+y+z$. The former gives the standard representation of $\mathrm{Sp}_4$ with trivial $\mathrm{SL}_2$-action, and the latter gives the external tensor product of the five-dimensional irreducible representation $W$ of $\mathrm{Sp}_4$ and the standard representation of $\mathrm{SL}_2$. This gives case (iii).
			\item At least two coordinate planes contain diagonal roots. After permuting the coordinates and changing their signs, we may assume that $x-y$ and $y-z$ are roots. The Weyl reflection $s_{x-y}$ sends $y-z$ to $x-z$. Thus all $\pm(x-y), \pm(y-z), \pm(x-z)$ are roots. They span a two-dimensional space. Since the rank of $\mathbf{G}_{\lambda}^{\mathrm{der}}$ is three, there is another root outside this span. If this additional root is diagonal, up to a sign it must be one of $x+y, y+z, z+x$. Hence by the claim above, the coordinate candidates in that plane are roots. So in either case we always have a coordinate root. Since the reflections $s_{x-y}$ and $s_{y-z}$ permute the three coordinate axes transitively, all six coordinate candidates are roots. Finally, the reflections $s_{2x}, s_{2y}, s_{2z}$ give the other diagonal line in each coordinate plane. Hence all diagonal candidates are roots. So all 18 candidates are roots, giving $C_3$ root system. Choose simple roots as $x-y,y-z,2z$. The highest weight in the formal character is $x+y+z=\omega_3$. This gives case (i).
		\end{itemize}
	\end{proof}
	
	In cases (ii) and (iii), for every irreducible constituent $U$ of $\rho_{\lambda}$, every irreducible summand of $U^{\mathrm{der}}$ is self-dual. Moreover, the derived monodromy group at $\lambda_0$ is of type $C_3$, whose Dynkin diagram has no nontrivial automorphisms, so $\mathbf{G}_{\lambda_0}^{\mathrm{der}}$ has no outer automorphisms. Therefore, condition $(\ast)$ of \autoref{Carbon} is satisfied. Since the derived formal character of $\rho_{\lambda_0}$ is multiplicity-free and $\rho_{\lambda}$ is totally odd, \autoref{Carbon}(ii) implies that every irreducible constituent of $\rho_{\lambda}$ is polarized and totally odd. Thus by \autoref{Magnesium}, the set of rational primes $\ell$ for which case (ii) or (iii) occurs for at least one $\lambda|\ell$ has Dirichlet density zero. This proves \autoref{Silicon}(ii).
	
	\subsection{Exterior powers of type $A$}
	
	In this subsection, we prove \autoref{Silicon}(iii). In fact, in this case, the irreducibility of every $\rho_{\lambda}$ follows directly from the geometry of the derived formal character, using the invariance-of-roots result of Hui \cite[Corollary 2.5]{Hu25}. We formulate the general result in the following lemma.
	\begin{lem}\label{Sulfur}
		Let $\mathbf{G}$ be a connected semisimple group over an algebraically closed field of characteristic zero. If a faithful representation $V$ of $G$ has the same formal character as $\bigwedge^k\mathrm{Std}$ of $\mathrm{SL}_n$, where $1\le k<n$, then
		$$(G,V)\simeq\left(\mathrm{SL}_n/\mu_{\gcd(n,k)},\bigwedge^k\mathrm{Std}\right)$$
		except possibly when $(n,k)=(4,2)$. In particular, if $(n,k)\not=(4,2)$, then $(G,V)$ is irreducible.
	\end{lem}
	\begin{proof}
		We identify $V$ with the underlying vector space of $\bigwedge^k\mathrm{Std}$ since they have the same dimension. Let $\mathbf{G}_0\subseteq GL_V$ be the image of $\mathrm{SL}_n$ under $\bigwedge^k\mathrm{Std}$.  Since $\mathbf{G}$ and $\mathbf{G}_0$ have the same formal character, after conjugation we may assume that $\mathbf{G}$ and $\mathbf{G}_0$ have a common maximal torus $T$. Identify $X^*(T)\otimes_{\mathbb{Z}}\mathbb{R}$ with the hyperplane
		$$H=\{(x_1,\cdots,x_n)\in\mathbb{R}^n:\sum_ix_i=0\},$$
		so that the root system of $\mathbf{G}_0$ is
		$$R_0=\{e_i-e_j:i\not=j\}.$$
		By \cite[Corollary 2.5]{Hu25}, $R_0$ is invariant under the action of $N_{\mathrm{GL}_V}(T)$ on $X^*(T)$, since $\mathbf{G}_0$ is almost simple of type $A_{n-1}$ and acts irreducibly on $V$. In particular, the Weyl group $W$ of $\mathbf{G}$ preserves $R_0$. Hence for $n\ge3$,
		$$W\subseteq\mathrm{Aut}(R_0)=\{\pm1\}\times S_n.$$
		We show that in fact $W\subseteq S_n$ for $n\ge5$: a coordinate permutation acts on $H$ as a reflection precisely when it is a transposition. If $-\sigma$ acts as a reflection on $H$, then $\sigma^2=\mathrm{id}$. Writing $\sigma$ as a product of $r$ disjoint transpositions, it follows that the dimension of the $(-1)$-eigenspace of $-\sigma$ on $H$ is $n-r-1$. Hence $r=n-2$, which is absurd for $n\ge5$, since $r\le\lfloor n/2\rfloor$.
		
		For $2\le n\le4$, our assumption implies that $\bigwedge^k\mathrm{Std}$ is either the standard representation of $\mathrm{SL}_n$ or its dual. In either case, the linear automorphisms of the weight set are precisely the coordinate permutations, so $W\subseteq S_n$.
		
		Hence if $(n,k)\not=(4,2)$, every Weyl reflection of $\mathbf{G}$ is a coordinate transposition. Since the roots span $H$, the graph of these transpositions on $\{1,\cdots,n\}$ is connected. Hence $W=S_n$. Thus every root of $\mathbf{G}$ lies on a line spanned by some $e_i-e_j$. Since every root is a difference of two weights of $V$, and the only nonzero weight differences on this line are $\pm(e_i-e_j)$, it follows that every root of $\mathbf{G}$ belongs to $R_0$. Since $W=S_n$ acts transitively on $R_0$ and permutes the roots of $\mathbf{G}$, the root system of $\mathbf{G}$ is precisely $R_0$. Since $\mathbf{G}$ and $\mathbf{G}_0$ have the same character lattice, formal character, roots and Weyl group, their root data and hence their tautological representations are isomorphic. This concludes the proof.
	\end{proof}
	
	By $\lambda$-independence of the derived formal character, \autoref{Sulfur} implies that the tautological representation of $\mathbf{G}_{\lambda}^{\mathrm{der}}$ is irreducible for every $\lambda$ in case \autoref{Silicon}(iii). Hence $\rho_{\lambda}$ is Lie-irreducible for every $\lambda$. This proves \autoref{Silicon}(iii).
	
	\section{Applications}
	
	\subsection{Large residual images}
	
	In this subsection, we prove \autoref{Argon}. Since twisting by a character does not change the projective image, we apply the normalization \autoref{Chlorine} to the cuspidal automorphic representation $\pi$. Then by our main theorem \autoref{main}, there is a Dirichlet density one set $\mathcal{L}$ of rational primes, such that $\mathbf{G}_{\lambda}\simeq G_2$ for every $\ell\in\mathcal{L}$ and every $\lambda|\ell$.
	
	We take $E=E_{\pi}$ to be the minimal coefficient field of the normalized compatible system. For a normalized automorphic compatible system of type $G_2$ with minimal coefficient field, we prove the following more precise large image result.
	\begin{prop}\label{Cobalt}
		Let $K$ be imaginary CM or totally real, and let $\pi$ be a regular algebraic polarized cuspidal automorphic representation of $\mathrm{GL}_7(\mathbb{A}_K)$ with associated $E$-rational compatible system $\{\rho_{\lambda}\}_{\lambda}$. Assume that for some coefficient prime $\lambda_0$, the algebraic monodromy group of $\rho_{\lambda_0}$ is connected semisimple of type $G_2$, and that $E$ is the minimal coefficient field of the system. Denote by $\kappa_{\lambda}$ the residue field of $E_{\lambda}$. Then there exists a Dirichlet density one set $\mathcal{L}$ of rational primes, such that for every $\ell\in\mathcal{L}$ and every $\lambda|\ell$,
		$$\overline{\rho}_{\lambda}^{\mathrm{ss}}(\mathrm{Gal}_K)=G_2(\kappa_{\lambda}),$$
		up to conjugacy, where $\overline{\rho}_{\lambda}^{\mathrm{ss}}$ denotes the semisimple reduction of $\rho_{\lambda}$.
	\end{prop}
	
	We first prove that, under our setting, the representation $\rho_{\lambda}$ factors through the coefficient field $E_{\lambda}$. The proof follows from the argument of Chenevier \cite[Corollary 6.7]{Ch19}.
	\begin{lem}
		For every coefficient prime $\lambda$ such that $\mathbf{G}_{\lambda}\simeq G_2$, the representation $\rho_{\lambda}$ can be defined over $E_{\lambda}$.
	\end{lem}
	\begin{proof}
		Fix a coefficient prime $\lambda$ such that $\mathbf{G}_{\lambda}\simeq G_2$. After conjugation, we may assume that $\rho_{\lambda}$ is valued in $G_2(\overline{E}_{\lambda})$ and has Zariski dense image. Fix $\sigma\in\mathrm{Gal}(\overline{E}_{\lambda}/E_{\lambda})$. Since the system is $E$-rational, Chebotarev and Brauer-Nesbitt imply that $^{\sigma}\rho_{\lambda}\simeq\rho_{\lambda}$. Hence there exists a matrix $A_{\sigma}\in\mathrm{GL}_7(\overline{E}_{\lambda})$ such that
		$$\sigma(\rho_{\lambda}(g))=A_{\sigma}^{-1}\rho_{\lambda}(g)A_{\sigma},\quad\forall g\in\mathrm{Gal}_K.$$
		Since $\mathrm{im}\rho_{\lambda}$ is Zariski-dense in $G_2/\overline{E}_{\lambda}$, $A_{\sigma}$ normalizes $G_2$. Since every automorphism of $G_2$ is inner and the standard representation of $G_2$ is irreducible, Schur's lemma allows us to replace $A_{\sigma}$ by an element $a_{\sigma}\in G_2(\overline{E}_{\lambda})$. Moreover, this element is unique since $G_2$ has trivial center. The uniqueness implies that for every $\sigma,\tau\in\mathrm{Gal}(\overline{E}_{\lambda}/E_{\lambda})$,
		$$a_{\sigma\tau}=a_{\sigma}\sigma(a_{\tau}).$$
		Hence $\sigma\mapsto a_{\sigma}$ is a Galois cocycle. It is moreover continuous for the discrete topology on $G_2(\overline{E}_{\lambda})$ since $\rho_{\lambda}$ is defined over a finite extension of $E_{\lambda}$. Since $G_2$ is simply connected and semisimple, \cite[Satz 1]{Kn65} gives $H^1(E_{\lambda},G_2)=1$. Hence there exists a matrix $A\in G_2(\overline{E}_{\lambda})$ such that $a_{\sigma}=A^{-1}\sigma(A)$ for every $\sigma\in\mathrm{Gal}(\overline{E}_{\lambda}/E_{\lambda})$. Therefore
		\begin{align*}
			\sigma(A\rho_{\lambda}(g)A^{-1})&=\sigma(A)\sigma(\rho_{\lambda}(g))\sigma(A)^{-1}\\&=Aa_{\sigma}\left(a_{\sigma}^{-1}\rho_{\lambda}(g)a_{\sigma}\right)a_{\sigma}^{-1}A^{-1}\\&=A\rho_{\lambda}(g)A^{-1}.
		\end{align*}
		Thus $A\rho_{\lambda}A^{-1}$ takes values in $G_2(E_{\lambda})$, as desired.
	\end{proof}
	
	To apply Larsen's maximality theorem, we first construct a $\mathbb{Q}$-rational compatible system by restriction of scalars. Put $d=[E:\mathbb{Q}]$.
	$$\left\{\rho_{\ell}=\bigoplus_{\lambda|\ell}\rho_{\lambda}:\mathrm{Gal}_K\to\left(\mathrm{Res}_{E/\mathbb{Q}}\mathrm{GL}_{7,E}\right)(\mathbb{Q}_{\ell})\subseteq\mathrm{GL}_{7d}(\mathbb{Q}_{\ell})\right\}_{\ell}.$$
	\begin{lem}
		Under our setting that $E$ is the minimal coefficient field, for every $\ell\in\mathcal{L}$, the algebraic monodromy group $\mathbf{H}_{\ell}$ of $\rho_{\ell}$ is
		$$\mathbf{H}_{\ell}=\prod_{\lambda|\ell}\mathrm{Res}_{E_{\lambda}/\mathbb{Q}_{\ell}}G_2.$$
	\end{lem}
	\begin{proof}
		By construction,
		\begin{equation}\label{7}
		\mathbf{H}_{\ell}\subseteq\prod_{\lambda|\ell}\mathrm{Res}_{E_{\lambda}/\mathbb{Q}_{\ell}}G_{2,E_{\lambda}}.
		\end{equation}
		Fix $\ell$ and put $X=\mathrm{Hom}_{\mathbb{Q}}(E,\overline{\mathbb{Q}}_{\ell})$. Extending the coefficient field to $\overline{\mathbb{Q}}_{\ell}$ yields
		$$\rho_{\ell}\otimes_{\mathbb{Q}_{\ell}}\overline{\mathbb{Q}}_{\ell}\simeq\bigoplus_{\sigma\in X}\rho_{\sigma},$$
		where $\rho_{\sigma}=\widetilde{\sigma}\circ\rho_{\lambda(\sigma)}$, where $\lambda(\sigma)$ is the coefficient prime induced by $\sigma$ and $\widetilde{\sigma}:E_{\lambda(\sigma)}\hookrightarrow\overline{\mathbb{Q}}_{\ell}$ the extension of $\sigma$. Each $\rho_{\sigma}$ has algebraic monodromy group $G_2/\overline{\mathbb{Q}}_{\ell}$. Hence
		$$\mathbf{H}_{\ell}\otimes_{\mathbb{Q}_{\ell}}\overline{\mathbb{Q}}_{\ell}\subseteq\prod_{\sigma\in X}G_2,$$
		which projects surjectively onto each factor.
		
		Suppose for a contradiction that this inclusion is proper. Since $G_2$ is simple and has trivial center, Goursat's lemma implies that there exist distinct $\sigma,\tau\in X$ and an algebraic automorphism $\alpha$ of $G_2$ such that $\rho_{\sigma}=\alpha\circ\rho_{\tau}$. Since every algebraic automorphism of $G_2$ is inner, we obtain $\rho_{\sigma}\simeq\rho_{\tau}$ as representations. In particular, for each finite place $v$ of $K$ outside the exceptional set and not above $\ell$,
		$$\sigma(P_v(T))=\tau(P_v(T)),$$
		where $P_v(T)$ is the common characteristic polynomial of $\rho_{\lambda}(\mathrm{Frob}_v)$ for almost all $\lambda$. Since $E$ is generated by coefficients of these $P_v(T)$, we have $\sigma=\tau$, a contradiction. Hence the inclusion (\ref{7}) is an equality after base change to $\overline{\mathbb{Q}}_{\ell}$, and therefore the inclusion is an equality.
	\end{proof}
	
	We recall Larsen's maximality theorem. Let 
	$$\rho_{\ell}:\mathrm{Gal}_K\to\mathrm{GL}_n(\mathbb{Q}_{\ell})$$
	be a semisimple continuous $\ell$-adic Galois representation of a number field $K$. We use the following notation:
	\begin{itemize}
		\item $\mathbf{H}_{\ell}$ the algebraic monodromy group of $\rho_{\ell}$,
		\item $\Gamma_{\ell}$ the Galois image $\rho_{\ell}(\mathrm{Gal}_K)$,
		\item $\mathbf{H}_{\ell}^{\mathrm{ad}}=\mathbf{H}_{\ell}^{\circ}/Z(\mathbf{H}_{\ell}^{\circ})$ the adjoint group of $\mathbf{H}_{\ell}^{\circ}$,
		\item $\mathbf{H}_{\ell}^{\mathrm{sc}}$ the simply connected cover of $\mathbf{H}_{\ell}^{\mathrm{ad}}$.
	\end{itemize}
	Consider the natural maps
	$$\xymatrix{\mathbf{H}_{\ell}^{\circ}\ar[r]^{q_{\ell}}&\mathbf{H}_{\ell}^{\mathrm{ad}}&\mathbf{H}_{\ell}^{\mathrm{sc}}\ar[l]_{p_{\ell}}},$$
	and define
	$$\Gamma_{\ell}^{\mathrm{sc}}=p_{\ell}^{-1}\left(q_{\ell}\left(\Gamma_{\ell}\cap\mathbf{H}_{\ell}^{\circ}(\mathbb{Q}_{\ell})\right)\right)\subseteq\mathbf{H}_{\ell}^{\mathrm{sc}}(\mathbb{Q}_{\ell}).$$
	\begin{thm}\cite[Theorem 3.17]{La95}
		Let $K$ be a number field and let
		$$\{\rho_{\ell}:\mathrm{Gal}_K\to\mathrm{GL}_n(\mathbb{Q}_{\ell})\}_{\ell\in\mathcal{L}}$$
		be a semisimple $\mathbb{Q}$-rational compatible system, where $\mathcal{L}$ is a set of rational primes of Dirichlet density one. Then there exists a density-one subset $\mathcal{L}'\subseteq\mathcal{L}$ such that $\Gamma_{\ell}^{\mathrm{sc}}$ is a hyperspecial maximal compact subgroup of $\mathbf{H}_{\ell}^{\mathrm{sc}}(\mathbb{Q}_{\ell})$ for every $\ell\in\mathcal{L}'$.
	\end{thm}
	
	Take the Dirichlet density one set $\mathcal{L}$ as in our setting. The theorem gives a density one subset $\mathcal{L}'\subseteq\mathcal{L}$. For each $\ell\in\mathcal{L}'$, the algebraic monodromy $\mathbf{H}_{\ell}$ is connected, simply connected with trivial center, hence $\mathbf{H}_{\ell}^{\mathrm{sc}}=\mathbf{H}_{\ell}$ and therefore $\Gamma_{\ell}^{\mathrm{sc}}=\Gamma_{\ell}$. Consequently, the Galois image $\Gamma_{\ell}$ itself is hyperspecial maximal compact subgroup of $\mathbf{H}_{\ell}(\mathbb{Q}_{\ell})$. By definition, let $\mathcal{H}'$ be a reductive $\mathbb{Z}_{\ell}$-model of $\mathbf{H}_{\ell}$ such that $\Gamma_{\ell}=\mathcal{H}'(\mathbb{Z}_{\ell})$. After removing finitely many primes from $\mathcal{L}'$, we may assume that $\ell$ is unramified in $E$. Let $\mathcal{H}$ be the reductive $\mathbb{Z}_{\ell}$-model
	$$\mathcal{H}=\prod_{\lambda|\ell}\mathrm{Res}_{\mathscr{O}_{E_{\lambda}}/\mathbb{Z}_{\ell}}G_2.$$
	Both models have generic fiber $\mathbf{H}_{\ell}$. Applying \cite[Theorem 7.2.16]{Co14} to the identity isomorphism of their generic fibers, and using $\mathbf{H}_{\ell}^{\mathrm{ad}}=\mathbf{H}_{\ell}$, we obtain an element $h\in\mathbf{H}_{\ell}(\mathbb{Q}_{\ell})$ such that conjugation by $h$ extends to an isomorphism $\mathcal{H}\simeq\mathcal{H}'$. Taking $\mathbb{Z}_{\ell}$-points yields
	$$\Gamma_{\ell}=h\mathcal{H}(\mathbb{Z}_{\ell})h^{-1}=h\left(\prod_{\lambda|\ell}G_2(\mathscr{O}_{E_{\lambda}})\right)h^{-1}.$$
	Write $h=(h_{\lambda})_{\lambda|\ell}$ under the identification
	$$\mathbf{H}_{\ell}(\mathbb{Q}_{\ell})=\prod_{\lambda|\ell}G_2(E_{\lambda}).$$
	Projecting onto the $\lambda$-factor, we obtain
	$$\rho_{\lambda}(\mathrm{Gal}_K)=h_{\lambda}G_2(\mathcal{O}_{E_{\lambda}})h_{\lambda}^{-1}.$$
	Hence the image of the semisimple reduction of $\rho_{\lambda}$ is $G_2(\kappa_{\lambda})$, as desired.
	
	Since enlarging the coefficient field to $\overline{\mathbb{Q}}_{\ell}$ does not change the representations, their semisimple reductions are isomorphic over $\overline{\mathbb{F}}_{\ell}$. Since twisting by a character does not change the residual projective image, \autoref{Cobalt} implies \autoref{Argon}.
	
	\subsection{Inverse Galois result}
	
	In this subsection, we prove \autoref{Calcium}. We first recall the following construction of Chenevier. For positive integers $u<v$, let $\mathcal{C}(u,v)$ denote the set of isomorphism classes of regular algebraic self-dual cuspidal automorphic representations $\pi$ of $\mathrm{GL}_7(\mathbb{A}_{\mathbb{Q}})$ of level one such that $\pi_{\infty}$ has weights $\{0,\pm u,\pm v,\pm(u+v)\}$.
	
	Let $\pi$ be a regular algebraic self-dual cuspidal automorphic representation of $\mathrm{GL}_7(\mathbb{A}_{\mathbb{Q}})$, and let $\{r_{\pi,\lambda}\}_{\lambda}$ be the associated compatible system as in \autoref{Titanium}. We set
	$$\rho_{\pi,\lambda}:=r_{\pi,\lambda}\otimes\varepsilon_{\ell}^3,$$
	where $\varepsilon_{\ell}$ is the $\ell$-adic cyclotomic character and $\lambda|\ell$. Then each $\rho_{\lambda}$ is self-dual. We call $\{\rho_{\pi,\lambda}\}_{\lambda}$ the normalized compatible system associated to $\pi$, as in \cite[\S 6.2]{Ch19}.
	\begin{thm}\cite[Theorems 6.4 and 6.12]{Ch19}\label{Chromium}
		For each pair
		$$(u,v)\in\{(4,8),(5,8),(3,10),(5,9),(4,10),(2,12),(7,8),(4,11),(1,14),(1,16),(1,17)\},$$
		we have $|\mathcal{C}(u,v)|=1$. Moreover, if $\pi$ is the unique representation in $\mathcal{C}(u,v)$, then for every rational prime $\ell$, its associated normalized Galois representation takes values, up to conjugacy, in the standard seven-dimensional representation of $G_2(\overline{\mathbb{Q}}_{\ell})$.
	\end{thm}
	
	Using Zenteno's results on large residual images \cite[Theorems 4.2 and 5.2, and Remark 2]{Zen23}, we show that, for each automorphic representation $\pi$ as above, the associated compatible system satisfies the hypothesis of our main theorem.
	\begin{prop}\label{Nickel}
		Let $\pi$ be one of the automorphic representations in \autoref{Chromium}. Then the associated normalized compatible system $\{\rho_{\pi,\ell}\}_{\ell}$ is $\mathbb{Q}$-rational, and for some rational prime $\ell_0$, the algebraic monodromy group $\mathbf{G}_{\ell_0}$ of $\rho_{\pi,\ell_0}$ is isomorphic to $G_2$.
	\end{prop}
	\begin{proof}
		By \cite[Theorem 3.13]{Cl90}, for every $\gamma\in\mathrm{Aut}(\mathbb{C})$, the Galois conjugate $^{\gamma}\pi$ is again a regular algebraic self-dual cuspidal automorphic representation of level one with the same algebraic weights as $\pi$. Hence $^{\gamma}\pi\in\mathcal{C}(u,v)$. Since $|\mathcal{C}(u,v)|=1$, it follows that $^{\gamma}\pi\simeq\pi$. Consequently, the coefficients of the Satake polynomials $\det(T-c(\pi_p))$ are fixed by every $\gamma\in\mathrm{Aut}(\mathbb{C})$, and hence belong to $\mathbb{Q}$. By \autoref{Titanium}(ii), these Satake polynomials are precisely the characteristic polynomials of the Frobenius elements under the normalized representation $\rho_{\pi,\ell}$. Hence the normalized system $\{\rho_{\pi,\ell}\}_{\ell}$ is $\mathbb{Q}$-rational.
		
		Next, we show that if $\pi$ is one of the eleven automorphic representations in \autoref{Chromium}, there exists a set $\mathcal{L}_{\pi}$ of rational primes of positive Dirichlet density, such that
		$$G_2(\mathbb{F}_{\ell})\subseteq\overline{\rho}_{\pi,\ell}^{\mathrm{ss}}(\mathrm{Gal}_{\mathbb{Q}})$$
		up to conjugacy for every $\ell\in\mathcal{L}_{\pi}$. For the ten pairs $(u,v)\not=(4,8)$ in \autoref{Chromium}, this follows from \cite[Theorem 4.2]{Zen23}, since $v\not=2u$. For the remaining pair $(4,8)$, this follows from \cite[Theorem 5.2 and Remark 2]{Zen23}.
		
		For a sufficiently large prime $\ell_0\in\mathcal{L}_{\pi}$, let $\underline{G}_{\ell_0}$ be the algebraic envelope of $\rho_{\pi,\ell_0}$ as in \autoref{Vanadium}, whose hypotheses are satisfied by \autoref{Iron}. Since the residual image contains $G_2(\mathbb{F}_{\ell_0})$, the Nori group contains $G_{2,\mathbb{F}_{\ell_0}}$. By \autoref{Vanadium}, after scalar extension if needed, the derived group of $\underline{G}_{\ell_0}$ therefore contains $G_{2,\mathbb{F}_{\ell_0}}$. Applying \autoref{Manganese} and using the inclusion $\mathbf{G}_{\ell}\subseteq G_2$ given by \autoref{Chromium}, we obtain
		$$14=\dim G_2\le\dim\underline{G}_{\ell_0}^{\mathrm{der}}\le\dim\underline{G}_{\ell_0}\le\dim\mathbf{G}_{\ell_0}^{\circ}\le\dim\mathbf{G}_{\ell_0}\le\dim G_2=14.$$
		Since $\mathbf{G}_{\ell_0}\subseteq G_2$ and $G_2$ is geometrically connected, the equality of dimensions implies $\mathbf{G}_{\ell_0}=G_2$, as desired. 
	\end{proof}
	
	Let $\pi_1,\cdots,\pi_{11}$ be the automorphic representations in \autoref{Chromium}. By \autoref{Cobalt} and \autoref{Nickel}, there exists a Dirichlet density one set $\mathcal{L}$ of rational primes such that, for every $\ell\in\mathcal{L}$ and every $1\le i\le 11$,
	$$\overline{\rho}_{\pi_i,\ell}^{\mathrm{ss}}(\mathrm{Gal}_{\mathbb{Q}})=G_2(\mathbb{F}_{\ell}),$$
	up to conjugacy in $\mathrm{GL}_7(\overline{\mathbb{F}}_{\ell})$. We may assume that $\ell>3$. As the representations $\pi_i$ have distinct algebraic weights, they are pairwise nonisomorphic. By strong multiplicity one, for each pair $i\not=j$, there exists a rational prime $p_{ij}$ at which their Satake polynomials are different. Since these polynomials have rational coefficients, after removing finitely many primes from $\mathcal{L}$, the residual representations $\overline{\rho}_{\pi_i,\ell}^{\mathrm{ss}}$ are pairwise nonisomorphic for every $\ell\in\mathcal{L}$.
	
	For each $\ell\in\mathcal{L}$, and each $1\le m\le11$, consider the product representation
	$$\overline{\rho}_{\ell,m}:=\prod_{1\le i\le m}\overline{\rho}_{\pi_i,\ell}^{\mathrm{ss}}:\mathrm{Gal}_{\mathbb{Q}}\to G_2(\mathbb{F}_{\ell})^m.$$
	Each projection of the image is surjective. Since $G_2(\mathbb{F}_{\ell})$ is nonabelian simple and every automorphism is inner for $\ell>3$ \cite[Table 5]{CCNPW85}, Goursat's lemma and the pairwise nonisomorphism of the residual representations imply that
	$$\overline{\rho}_{\ell,m}(\mathrm{Gal}_{\mathbb{Q}})=G_2(\mathbb{F}_{\ell})^m.$$
	Let $L_m$ be the fixed field of the kernel of $\overline{\rho}_{\ell,m}$. Then
	$$\mathrm{Gal}(L_m/\mathbb{Q})\simeq G_2(\mathbb{F}_{\ell})^m.$$
	Since all the $\pi_i$ are of level one, the representations $\rho_{\pi_i,\ell}$ are unramified outside $\ell$. Thus $L_m/\mathbb{Q}$ is unramified at every rational prime other than $\ell$, which proves \autoref{Calcium}.

	\section*{Acknowledgement}


\end{document}